\documentclass[a4paper, 12 pt]{article}

\usepackage[T2A]{fontenc}
\usepackage[cp1251]{inputenc}
\usepackage[english]{babel}
\usepackage[tbtags]{amsmath}
\usepackage{amsfonts,amssymb}

\usepackage{mathrsfs}

\usepackage{geometry} 
\begin{document}
\title{Riordan arrays and generalized Euler polynomials}
 \author{E. Burlachenko}
 \date{}

 \maketitle
\begin{abstract}
Generalization of the Euler polynomials ${{A}_{n}}\left( x \right)={{\left( 1-x \right)}^{n+1}}\sum\nolimits_{m=0}^{\infty }{{{m}^{n}}{{x}^{m}}}$ are the polynomials ${{\alpha }_{n}}\left( x \right)={{\left( 1-x \right)}^{n+1}}\sum\nolimits_{m=0}^{\infty }{{{u}_{n}}}\left( m \right){{x}^{m}}$, where ${{u}_{n}}\left( x \right)$ is the polynomial of degree $n$. These polynomials appear in various fields of mathematics, which causes a variety of methods for their study. In present  paper we will consider generalized Euler polynomials as an attribute of the theory of Riordan arrays. From this point of view, we will consider the transformations associated with them, with a participation of such objects as binomial sequences, Stirling numbers, multinomial coefficients, shift operator, and demonstrate a constructiveness of the chosen point of view.
\end{abstract}
\section{Introduction}

Transformations, corresponding to multiplication and composition of  series, play the main role in the space of formal power series over the field of real or complex  numbers. Multiplication is given by the matrix $\left( f\left( x \right),x \right)$ $n$th column of which, $n=0,\text{ }1,\text{ }2,\text{ }...$ ,  has the generating function $f\left( x \right){{x}^{n}}$; composition is given by the matrix $\left( 1,g\left( x \right) \right)$ $n$th column of which  has the generating function ${{g}^{n}}\left( x \right)$, ${{g}_{0}}=0$:
$$\left( f\left( x \right),x \right)a\left( x \right)=f\left( x \right)a\left( x \right),\qquad \left( 1,g\left( x \right) \right)a\left( x \right)=a\left( g\left( x \right) \right).$$
Matrix
$$\left( f\left( x \right),x \right)\left( 1,g\left( x \right) \right)=\left( f\left( x \right),g\left( x \right) \right)$$
is called Riordan array [1] – [5]; $n$th column of the Riordan array has the generating function $f\left( x \right){{g}^{n}}\left( x \right)$. Thus
$$\left( f\left( x \right),g\left( x \right) \right)b\left( x \right){{a}^{n}}\left( x \right)=f\left( x \right)b\left( g\left( x \right) \right){{a}^{n}}\left( g\left( x \right) \right),$$
$$\left( f\left( x \right),g\left( x \right) \right)\left( b\left( x \right),a\left( x \right) \right)=\left( f\left( x \right)b\left( g\left( x \right) \right),a\left( g\left( x \right) \right) \right).$$
Matrices $\left( f\left( x \right),g\left( x \right) \right)$, ${{f}_{0}}\ne 0$, ${{g}_{1}}\ne 0$, form a group, called the Riordan group. 

$n$th coefficient of the series $a\left( x \right)$, $n$th row and $n$th column of the matrix $A$ will be denoted  respectively by
$$\left[ {{x}^{n}} \right]a\left( x \right),   \qquad\left[ n,\to  \right]A,  \qquad[\uparrow ,n]A,$$
and
$$\left[ {{x}^{n}} \right]a\left( x \right)b\left( x \right)=\left[ {{x}^{n}} \right]\left( a\left( x \right)b\left( x \right) \right).$$
We associate rows and columns of matrices with the generating functions of their elements. Matrices 
$${{\left| {{e}^{x}} \right|}^{-1}}\left( f\left( x \right),g\left( x \right) \right)\left| {{e}^{x}} \right|={{\left( f\left( x \right),g\left( x \right) \right)}_{{{e}^{x}}}},$$
where $\left| {{e}^{x}} \right|$ is the diagonal matrix whose diagonal elements are equal to the coefficients of  the series ${{e}^{x}}$: $\left| {{e}^{x}} \right|a\left( x \right)=\sum\nolimits_{n=0}^{\infty }{{{{a}_{n}}{{x}^{n}}}/{n!}\;}$, are called exponential Riordan arrays. Denote
$$\left[ n,\to  \right]{{\left( f\left( x \right),g\left( x \right) \right)}_{{{e}^{x}}}}={{s}_{n}}\left( x \right), \qquad{{f}_{0}}\ne 0, \qquad{{g}_{1}}\ne 0.$$
Then
$${{\left( f\left( x \right),g\left( x \right) \right)}_{{{e}^{x}}}}{{\left( 1-\varphi x \right)}^{-1}}={{\left| {{e}^{x}} \right|}^{-1}}\left( f\left( x \right),g\left( x \right) \right){{e}^{\varphi x}}={{\left| {{e}^{x}} \right|}^{-1}}f\left( x \right)\exp \left( \varphi g\left( x \right) \right)$$
or
$$\sum\limits_{n=0}^{\infty }{\frac{{{s}_{n}}\left( \varphi  \right)}{n!}{{x}^{n}}}=f\left( x \right)\exp \left( \varphi g\left( x \right) \right).$$
Sequence of polynomials ${{s}_{n}}\left( x \right)$ is called Sheffer sequence, and in the case $f\left( x \right)=1$, binomial sequence. Matrix
$$P=\left( \frac{1}{1-x},\frac{x}{1-x} \right)={{\left( {{e}^{x}},x \right)}_{{{e}^{x}}}}=\left( \begin{matrix}
   1 & 0 & 0 & 0 & 0 & \ldots   \\
   1 & 1 & 0 & 0 & 0 & \ldots   \\
   1 & 2 & 1 & 0 & 0 & \ldots   \\
   1 & 3 & 3 & 1 & 0 & \ldots   \\
   1 & 4 & 6 & 4 & 1 & \ldots   \\
   \vdots  & \vdots  & \vdots  & \vdots  & \vdots  & \ddots   \\
\end{matrix} \right)$$
is called Pascal matrix. Power of the Pascal matrix is defined by the identity
$${{P}^{\varphi }}=\left( \frac{1}{1-\varphi x},\frac{x}{1-\varphi x} \right)={{\left( {{e}^{\varphi x}},x \right)}_{{{e}^{x}}}}.$$

Paper [1]  laid the foundations for the theory of Riordan arrays and introduced a terminology that became generally accepted.  Prior to this, Riordan matrices and  their varieties were considered in the literature under different names. Series of papers [6] – [16] is devoted to study of matrices called convolution arrays. $n$th column of the convolution array has the generating function $b\left( x \right){{a}^{n}}\left( x \right)$, where ${{a}_{0}}=1$ or ${{a}_{0}}=0$, depending on the problem under consideration. Results obtained for the convolution arrays in terms of the Riordan arrays can be stated more concisely, but for this the constraint ${{a}_{0}}=0$ for the matrix $\left( b\left( x \right),a\left( x \right) \right)$ must be removed. Thus, along with lower triangular  Riordan  matrices, we will consider “square” matrices $\left( b\left( x \right),a\left( x \right) \right)$, ${{a}_{0}}=1$. For example,
$$\left( 1,\frac{1}{1+x} \right)=\left( \begin{matrix}
   1 & 1 & 1 & 1 & \cdots   \\
   0 & -1 & -2 & -3 & \cdots   \\
   0 & 1 & 3 & 6 & \cdots   \\
   0 & -1 & -4 & -10 & \cdots   \\
   \vdots  & \vdots  & \vdots  & \vdots  & \ddots   \\
\end{matrix} \right).$$
This includes the upper triangular matrix $\left( 1,1+x \right)$, whose transpose is the Pascal matrix and which coincides with the  matrix of shift operator:
$$\left( 1,1+x \right)={{P}^{T}}=E=\left( \begin{matrix}
   1 & 1 & 1 & 1 & \cdots   \\
   0 & 1 & 2 & 3 & \cdots   \\
   0 & 0 & 1 & 3 & \cdots   \\
   0 & 0 & 0 & 1 & \cdots   \\
   \vdots  & \vdots  & \vdots  & \vdots  & \ddots   \\
\end{matrix} \right).$$
Matrix $\left( b\left( x \right),a\left( x \right) \right)$ can be multiplied from the right by the matrix with the finite columns and from the left by the matrix with the finite rows. If $c\left( x \right)$, $d\left( x \right)$ are polynomials, then
$$\left( b\left( x \right),a\left( x \right) \right)\left( c\left( x \right),d\left( x \right) \right)=\left( b\left( x \right)c\left( a\left( x \right) \right),d\left( a\left( x \right) \right) \right);$$
if ${{g}_{0}}=0$, then
$$\left( f\left( x \right),g\left( x \right) \right)\left( b\left( x \right),a\left( x \right) \right)=\left( f\left( x \right)b\left( g\left( x \right) \right),a\left( g\left( x \right) \right) \right)$$
(see [17], where square Riordan arrays are called generalized Riordan arrays). Denote
$$\left[ n,\to  \right]\left( 1,a\left( x \right)-1 \right)={{v}_{n}}\left( x \right)=\sum\limits_{m=1}^{n}{{{v}_{m}}{{x}^{m}}}, \qquad n>0.$$
Identities
$$\left( 1,a\left( x \right)-1 \right)\left( 1,1+x \right)=\left( 1,a\left( x \right) \right),  \qquad\left( 1,a\left( x \right)-1 \right)\left( 1,\frac{1}{1+x} \right)=\left( 1,{{a}^{-1}}\left( x \right) \right)$$
bear the following information. Since
$$\left[ n,\to  \right]\left( 1,1+x \right)=\frac{{{x}^{n}}}{{{\left( 1-x \right)}^{n+1}}};  \qquad\left[ n,\to  \right]\left( 1,\frac{1}{1+x} \right)=\frac{{{\left( -1 \right)}^{n}}x}{{{\left( 1-x \right)}^{n+1}}}, \qquad n>0,$$
then
$$\left[ n,\to  \right]\left( 1,a\left( x \right) \right)=\frac{{{\alpha }_{n}}\left( x \right)}{{{\left( 1-x \right)}^{n+1}}},  \qquad\left[ n,\to  \right]\left( 1,{{a}^{-1}}\left( x \right) \right)=\frac{\alpha _{n}^{\left( -1 \right)}\left( x \right)}{{{\left( 1-x \right)}^{n+1}}},$$
$${{\alpha }_{n}}\left( x \right)=\sum\limits_{m=1}^{n}{{{v}_{m}}{{x}^{m}}{{\left( 1-x \right)}^{n-m}}},  \qquad\alpha _{n}^{\left( -1 \right)}\left( x \right)=\sum\limits_{m=1}^{n}{{{v}_{m}}}{{\left( -1 \right)}^{m}}x{{\left( 1-x \right)}^{n-m}},$$
whence 
	$$\alpha _{n}^{\left( -1 \right)}\left( x \right)={{\left( -1 \right)}^{n}}x{{\hat{I}}_{n}}{{\alpha }_{n}}\left( x \right),\eqno	(1)$$	
where ${{\hat{I}}_{n}}$ is the operator (matrix) exchanging the coefficients of the polynomial of degree $n$ in reverse order. For example,
$${{\hat{I}}_{3}}=\left( \begin{matrix}
   0 & 0 & 0 & 1  \\
   0 & 0 & 1 & 0  \\
   0 & 1 & 0 & 0  \\
   1 & 0 & 0 & 0  \\
\end{matrix} \right).$$

If $a\left( x \right)={{e}^{x}}$, then ${{\alpha }_{n}}\left( x \right)={{{A}_{n}}\left( x \right)}/{n!}\;$, where ${{A}_{n}}\left( x \right)$ are the Euler polynomials:
$$\frac{{{A}_{n}}\left( x \right)}{{{\left( 1-x \right)}^{n+1}}}=\sum\limits_{m=0}^{\infty }{{{m}^{n}}}{{x}^{m}},  \qquad{{A}_{n}}\left( 1 \right)=n!.$$
For example,
$${{A}_{1}}\left( x \right)=x,  \qquad{{A}_{2}}\left( x \right)=x+{{x}^{2}},  \qquad{{A}_{3}}\left( x \right)=x+4{{x}^{2}}+{{x}^{3}},$$
$${{A}_{4}}\left( x \right)=x+11{{x}^{2}}+11{{x}^{3}}+{{x}^{4}}.$$

Polynomials associated with the generating functions of the rows of the convolution arrays (they are called numerator polynomials) are considered in papers [7], [8], [13], [14], [16]. Focus is not on general issues (except paper [8]), but on specific cases associated with the Fibonacci, Catalan sequences and their generalizations. Two such examples will be considered in this paper (Example 3, Example 7).

Concept of the generalized Euler polynomials in general form is represented in [18]. Polynomials under consideration (they are denoted as well as ordinary Euler polynomials) are called ${{p}_{n}}$-associated Eulerian polynomials and are defined as follows:
$$\frac{{{A}_{n}}\left( x \right)}{{{\left( 1-x \right)}^{n+1}}}=\sum\limits_{m=0}^{\infty }{{{p}_{n}}}\left( m \right){{x}^{m}},  \qquad{{p}_{n}}\left( x \right)=\sum\limits_{m=0}^{n}{\left( \begin{matrix}
   x+n-m  \\
   n  \\
\end{matrix} \right)}\left[ {{x}^{m}} \right]{{A}_{n}}\left( x \right).$$
Any polynomial sequence can be taken as sequence of the polynomials ${{p}_{n}}\left( x \right)$, but the most interesting case arises when ${{p}_{n}}\left( x \right)$ is a Sheffer sequence. In this case
$$\sum\limits_{n=0}^{\infty }{{{p}_{n}}}\left( t \right)\frac{{{x}^{n}}}{n!}=g\left( x \right)\exp \left( tf\left( x \right) \right),  \quad{{g}_{0}}\ne 0, \quad{{f}_{0}}=0, \quad{{f}_{1}}\ne 0,$$
$$\sum\limits_{n=0}^{\infty }{{{A}_{n}}}\left( t \right)\frac{{{x}^{n}}}{n!}=g\left( \left( 1-t \right)x \right)\frac{1-t}{1-t\exp \left( f\left( \left( 1-t \right)x \right) \right)}.\eqno 	(2)$$
In terms of the Riordan arrays this means that
$${{p}_{n}}\left( x \right)=\left[ n,\to  \right]{{\left( g\left( x \right),f\left( x \right) \right)}_{{{e}^{x}}}},    \qquad\frac{{{A}_{n}}\left( x \right)}{n!{{\left( 1-x \right)}^{n+1}}}=\left[ n,\to  \right]\left( g\left( x \right),{{e}^{f\left( x \right)}} \right).$$
We narrow the scope of this generalization and will consider generalized Euler polynomials (GEP) as ${{p}_{n}}$-associated Eulerian polynomials when ${{p}_{n}}\left( x \right)$ is a binomial sequence.

In Section 2 we consider the basic transformations associated with the GEP. Denote (we will bear in mind that always $n>0$):
$$\left[ n,\to  \right]\left( 1,a\left( x \right) \right)=\frac{{{\alpha }_{n}}\left( x \right)}{{{\left( 1-x \right)}^{n+1}}},  \qquad\left[ n,\to  \right]{{\left( 1,\log a\left( x \right) \right)}_{{{e}^{x}}}}={{u}_{n}}\left( x \right),$$
$$\left[ n,\to  \right]\left( 1,a\left( x \right)-1 \right)={{v}_{n}}\left( x \right),$$
$$\frac{1}{x}{{\alpha }_{n}}\left( x \right)={{\tilde{\alpha }}_{n}}\left( x \right),  \qquad\frac{1}{x}{{A}_{n}}\left( x \right)={{\tilde{A}}_{n}}\left( x \right),  \qquad\frac{1}{x}{{u}_{n}}\left( x \right)={{\tilde{u}}_{n}}\left( x \right),$$ 
$$\frac{1}{x}{{v}_{n}}\left( x \right)={{\tilde{v}}_{n}}\left( x \right),   \qquad{{\tilde{I}}_{n}}={{\hat{I}}_{n-1}}.$$
We introduce the matrices ${{U}_{n}}$, ${{V}_{n}}$, ${{V}_{n}}{{U}_{n}}$:
$$[\uparrow ,p]{{U}_{n}}=\frac{1}{n!}{{\left( 1-x \right)}^{n-1-p}}{{\tilde{A}}_{p+1}}\left( x \right),  \qquad[\uparrow ,p]U_{n}^{-1}=\frac{1}{x}\prod\limits_{m=0}^{n-1}{\left( x-p+m \right)},$$
$$[\uparrow ,p]{{V}_{n}}={{\left( 1+x \right)}^{n-p-1}}{{x}^{p}},    \qquad[\uparrow ,p]V_{n}^{-1}={{\left( 1-x \right)}^{n-p-1}}{{x}^{p}},$$
$$[\uparrow ,p]\left( {{V}_{n}}{{U}_{n}} \right)=\frac{1}{n!}\sum\limits_{m=1}^{p+1}{m!S\left( p+1,\text{ }m \right)}\text{ }{{x}^{m-1}},$$
$$[\uparrow ,p]\left( U_{n}^{-1}V_{n}^{-1} \right)=\frac{n!}{\left( p+1 \right)!}\sum\limits_{m=1}^{p+1}{s\left( p+1,\text{ }m \right){{x}^{m-1}}},$$
$p=0$, $1$, … , $n-1$; $S\left( p+1,\text{ }m \right)$ are the Stirling numbers of the second kind, $s\left( p+1,\text{ }m \right)$ are the Stirling numbers of  the first kind. Then
$${{U}_{n}}{{\tilde{u}}_{n}}\left( x \right)={{\tilde{\alpha }}_{n}}\left( x \right),  \qquad{{V}_{n}}{{\tilde{\alpha }}_{n}}\left( x \right)={{V}_{n}}{{U}_{n}}{{\tilde{u}}_{n}}\left( x \right)={{\tilde{v}}_{n}}\left( x \right).$$

In Section 3 we consider examples of an application of these transformations; in Example 4 we introduce analog of the GEP for a formal Dirichlet series.

In Section 4 we consider the transformation
$${{W}_{\left( n,\text{ }m \right)}}={{U}_{n}}\left( m,mx \right)U_{n}^{-1}.$$
Elements of the matrix ${{W}_{\left( n,\text{ }m \right)}}$ are the multinomial coefficients. Namely, let ${{\left( b\left( x \right),x \right)}_{m}}$ is the matrix such that
$$\left[ n,\to  \right]{{\left( b\left( x \right),x \right)}_{m}}=\left[ mn+m-1,\to  \right]\left( b\left( x \right),x \right).$$
For example,
$${{\left( b\left( x \right),x \right)}_{2}}=\left( \begin{matrix}
   {{b}_{1}} & {{b}_{0}} & 0 & 0 & \cdots   \\
   {{b}_{3}} & {{b}_{2}} & {{b}_{1}} & {{b}_{0}} & \cdots   \\
   {{b}_{5}} & {{b}_{4}} & {{b}_{3}} & {{b}_{2}} & \cdots   \\
   {{b}_{7}} & {{b}_{6}} & {{b}_{5}} & {{b}_{4}} & \cdots   \\
   \vdots  & \vdots  & \vdots  & \vdots  & \ddots   \\
\end{matrix} \right),  \qquad{{\left( b\left( x \right),x \right)}_{3}}=\left( \begin{matrix}
   {{b}_{2}} & {{b}_{1}} & {{b}_{0}} & 0 & \cdots   \\
   {{b}_{5}} & {{b}_{4}} & {{b}_{3}} & {{b}_{2}} & \cdots   \\
   {{b}_{8}} & {{b}_{7}} & {{b}_{6}} & {{b}_{5}} & \cdots   \\
   {{b}_{11}} & {{b}_{10}} & {{b}_{9}} & {{b}_{8}} & \cdots   \\
   \vdots  & \vdots  & \vdots  & \vdots  & \ddots   \\
\end{matrix} \right).$$
Then
$${{W}_{\left( n,\text{ }m \right)}}={{\left( {{\left( \frac{1-{{x}^{m}}}{1-x} \right)}^{n+1}},x \right)}_{m}}{{I}_{n}},$$
where ${{I}_{n}}$ is the Identity square matrix of order $n$. Matrices $\left( {1}/{{{m}^{n}}}\; \right){{\left( {{W}_{\left( n,m \right)}} \right)}^{T}}$ are known as “amazing matrices” [19, pp. 156-160], [20] – [23]. They find application in various fields of mathematics. From point of view of the theory of Riordan arrays, the transformation ${{W}_{\left( n,m \right)}}$ has the following sense. Denote
$$\left[ n,\to  \right]\left( 1,{{a}^{m}}\left( x \right) \right)=\frac{\alpha _{n}^{\left( m \right)}\left( x \right)}{{{\left( 1-x \right)}^{n+1}}},  \qquad\alpha _{n}^{\left( 1 \right)}\left( x \right)=\alpha _{n}^{{}}\left( x \right),  \qquad\frac{1}{x}\alpha _{n}^{\left( m \right)}\left( x \right)=\tilde{\alpha }_{n}^{\left( m \right)}\left( x \right).$$
Then
$${{W}_{\left( n,m \right)}}{{\tilde{\alpha }}_{n}}\left( x \right)=\tilde{\alpha }_{n}^{\left( m \right)}\left( x \right).$$

In Section 5 we consider the transformation $A_{n}^{\beta }={{U}_{n}}{{E}^{n\beta }}U_{n}^{-1}$, which has the following sense. Each formal power series $a\left( x \right)$, ${{a}_{0}}=1$, is associated by means of the Lagrange transform
$${{a}^{\varphi }}\left( x \right)=\sum\limits_{n=0}^{\infty }{\frac{{{x}^{n}}}{{{a}^{\beta n}}\left( x \right)}\left[ {{x}^{n}} \right]}\left( 1-x\beta {{\left( \log a\left( x \right) \right)}^{\prime }} \right){{a}^{\varphi +\beta n}}\left( x \right)$$
with the set of series $_{\left( \beta  \right)}a\left( x \right)$, $_{\left( 0 \right)}a\left( x \right)=a\left( x \right)$,   such that
$${}_{\left( \beta  \right)}a\left( x{{a}^{-\beta }}\left( x \right) \right)=a\left( x \right),   \qquad a\left( x{}_{\left( \beta  \right)}{{a}^{\beta }}\left( x \right) \right)={}_{\left( \beta  \right)}a\left( x \right),$$
$$\left[ {{x}^{n}} \right]{}_{\left( \beta  \right)}{{a}^{\varphi }}\left( x \right)=\left[ {{x}^{n}} \right]\left( 1-x\beta \frac{{a}'\left( x \right)}{a\left( x \right)} \right){{a}^{\varphi +\beta n}}\left( x \right)=\frac{\varphi }{\varphi +\beta n}\left[ {{x}^{n}} \right]{{a}^{\varphi +\beta n}}\left( x \right),$$
$$\left[ {{x}^{n}} \right]\left( 1+x\beta \frac{_{\left( \beta  \right)}{a}'\left( x \right)}{_{\left( \beta  \right)}a\left( x \right)} \right){}_{\left( \beta  \right)}{{a}^{\varphi }}\left( x \right)=\frac{\varphi +\beta n}{\varphi }\left[ {{x}^{n}} \right]{}_{\left( \beta  \right)}{{a}^{\varphi }}\left( x \right)=\left[ {{x}^{n}} \right]a_{{}}^{\varphi +\beta n}\left( x \right).$$
Denote
$$\left[ n,\to  \right]\left( 1,{}_{\left( \beta  \right)}a\left( x \right) \right)=\frac{_{\left( \beta  \right)}{{\alpha }_{n}}\left( x \right)}{{{\left( 1-x \right)}^{n+1}}},  \qquad\frac{1}{x}{}_{\left( \beta  \right)}{{\alpha }_{n}}\left( x \right)={}_{\left( \beta  \right)}{{\tilde{\alpha }}_{n}}\left( x \right).$$
Then
$$A_{n}^{\beta }{{\tilde{\alpha }}_{n}}\left( x \right)={}_{\left( \beta  \right)}{{\tilde{\alpha }}_{n}}\left( x \right).$$
It is interesting to observe how the properties of the shift operator, inherited by the transformation $A_{n}^{\beta }$, manifest themselves in new qualities. For example, since 
 $${{U}_{n}}\left( 1,-x \right)U_{n}^{-1}={{\left( -1 \right)}^{n+1}}{{\tilde{I}}_{n}},$$
then the transformations ${{\tilde{I}}_{n}}A_{n}^{\beta }$, $A_{n}^{\beta }{{\tilde{I}}_{n}}$ are involutions and ${{\tilde{I}}_{n}}A_{n}^{\beta }{{\tilde{I}}_{n}}=A_{n}^{-\beta }$. In Example 7, following [13], we give a general formula for the GEP associated with the generalized binomial series. Namely, let
$$_{\left( \beta  \right)}{{a}^{m}}\left( x \right)=\sum\limits_{n=0}^{\infty }{\frac{m}{m+\beta n}}\left( \begin{matrix}
   m+\beta n  \\
   n  \\
\end{matrix} \right){{x}^{n}},$$
$$\frac{_{\left( \beta  \right)}{{\alpha }_{n}}\left( x \right)}{{{\left( 1-x \right)}^{n+1}}}=\sum\limits_{m=0}^{\infty }{\frac{m}{m+\beta n}}\left( \begin{matrix}
   m+\beta n  \\
   n  \\
\end{matrix} \right){{x}^{m}}.$$
Then 
$$_{\left( \beta  \right)}{{\alpha }_{n}}\left( x \right)=\frac{1}{n}\sum\limits_{m=1}^{n}{\left( \begin{matrix}
   n\left( 1-\beta  \right)  \\
   m-1  \\
\end{matrix} \right)\left( \begin{matrix}
   n\beta   \\
   n-m  \\
\end{matrix} \right){{x}^{m}}}.$$
In particular,
$$_{\left( 0 \right)}{{\alpha }_{n}}\left( x \right)={{x}^{n}}, \qquad_{\left( 1 \right)}{{\alpha }_{n}}\left( x \right)=x,  \qquad_{\left( {1}/{2}\; \right)}{{\alpha }_{2n}}\left( x \right)=\frac{1}{2}\left( 1+x \right){{x}^{n}}.$$
\section{Basic transformations}

Denote
$$\left[ n,\to  \right]{{\left( 1,\log a\left( x \right) \right)}_{{{e}^{x}}}}={{u}_{n}}\left( x \right)=\sum\limits_{p=1}^{n}{{{u}_{p}}}{{x}^{p}}.$$
Then
$${{a}^{m}}\left( x \right)=\sum\limits_{n=0}^{\infty }{\frac{{{u}_{n}}\left( m \right)}{n!}{{x}^{n}}},   \qquad\frac{{{\alpha }_{n}}\left( x \right)}{{{\left( 1-x \right)}^{n+1}}}=\sum\limits_{m=0}^{\infty }{\frac{{{u}_{n}}\left( m \right)}{n!}}{{x}^{m}},$$
$$\frac{{{\alpha }_{n}}\left( x \right)}{{{\left( 1-x \right)}^{n+1}}}=\frac{1}{n!}\sum\limits_{m=0}^{\infty }{{{x}^{m}}}\sum\limits_{p=1}^{n}{{{u}_{p}}}{{m}^{p}}=\frac{1}{n!}\sum\limits_{p=1}^{n}{\sum\limits_{m=0}^{\infty }{{{u}_{p}}{{m}^{p}}{{x}^{m}}}}=$$
$$=\frac{1}{n!}\sum\limits_{p=1}^{n}{\frac{{{u}_{p}}{{A}_{p}}\left( x \right)}{{{\left( 1-x \right)}^{p+1}}}}=\frac{\frac{1}{n!}\sum\limits_{p=1}^{n}{{{u}_{p}}{{\left( 1-x \right)}^{n-p}}{{A}_{p}}\left( x \right)}}{{{\left( 1-x \right)}^{n+1}}}.$$
We introduce the matrices ${{U}_{n}}$:
$$[\uparrow ,p]{{U}_{n}}=\frac{1}{n!}{{\left( 1-x \right)}^{n-1-p}}{{\tilde{A}}_{p+1}}\left( x \right),  \qquad p=0, 1, … , n-1.$$
For example,
$$U_{2}^{{}}=\frac{1}{2}\left( \begin{matrix}
   1 & 1  \\
   -1 & 1  \\
\end{matrix} \right), \qquad U_{3}^{{}}=\frac{1}{3!}\left( \begin{matrix}
  1 & 1 & 1  \\
   -2 & 0 & 4  \\
  1 & -1 & 1  \\
\end{matrix} \right), \qquad U_{4}^{{}}=\frac{1}{4!}\left( \begin{matrix}
   1 & 1 & 1 & 1  \\
   -3 & -1 & 3 & 11  \\
   3 & -1 & -3 & 11  \\
   -1 & 1 & -1 & 1  \\
\end{matrix} \right).$$
Then
$${{U}_{n}}{{\tilde{u}}_{n}}\left( x \right)={{\tilde{\alpha }}_{n}}\left( x \right).$$
{\bfseries Theorem 1.}
$${{U}_{n}}\left( 1,-x \right)={{\left( -1 \right)}^{n+1}}{{\tilde{I}}_{n}}{{U}_{n}}.$$
{\bfseries Proof.} This follows from the identities (1) and 
$$\left[ n,\to  \right]{{\left( 1,\log {{a}^{-1}}\left( x \right) \right)}_{{{e}^{x}}}}={{u}_{n}}\left( -x \right).$$
{\bfseries Theorem 2.}
$${{\alpha }_{n}}\left( 1 \right)={{\left( {{a}_{1}} \right)}^{n}}.$$
{\bfseries Proof.} Denote $[\uparrow ,p]{{U}_{n}}={{U}_{p}}\left( x \right)$. Since
$${{a}_{1}}=\left[ x \right]\log a\left( x \right), \qquad{{\left( {{a}_{1}} \right)}^{n}}=\left[ {{x}^{n}} \right]{{u}_{n}}\left( x \right);  \qquad{{U}_{p}}\left( 1 \right)=0, \qquad p<n-1;   \qquad{{U}_{n-1}}\left( 1 \right)=1.$$
then
$${{\alpha }_{n}}\left( 1 \right)=\sum\limits_{p=0}^{n-1}{{{u}_{p+1}}{{U}_{p}}}\left( 1 \right)={{u}_{n}}.$$ 
{\bfseries Theorem 3.}
$$[\uparrow ,p]U_{n}^{-1}=\frac{1}{x}\prod\limits_{m=0}^{n-1}{\left( x-p+m \right)}, \qquad p=0, 1, … , n-1.$$
{\bfseries Proof.} Denote 
$$\frac{{{x}^{p}}}{{{\left( 1-x \right)}^{n+1}}}=\sum\limits_{m=0}^{\infty }{\frac{^{\left( p \right)}{{u}_{n}}\left( m \right)}{n!}{{x}^{m}}},  \qquad p=1, 2, …, n.$$
Then, according to (1),
$$\frac{{{\left( -1 \right)}^{n}}x{{{\hat{I}}}_{n}}{{x}^{p}}}{{{\left( 1-x \right)}^{n+1}}}=\frac{{{\left( -1 \right)}^{n}}{{x}^{n-p+1}}}{{{\left( 1-x \right)}^{n+1}}}=\sum\limits_{m=0}^{\infty }{\frac{^{\left( p \right)}{{u}_{n}}\left( -m \right)}{n!}{{x}^{m}}}.$$
Hence, $^{\left( p \right)}{{u}_{n}}\left( x \right)=0$ if $x=p-1$, $p-2$, …, $p-n$. I.e.  ${}^{\left( p \right)}{{u}_{n}}\left( x \right)=\prod\nolimits_{m=1}^{n}{\left( x-p+m \right)}$. Let ${{\alpha }_{n}}\left( x \right)=\sum\nolimits_{p=1}^{n}{{{\alpha }_{p}}{{x}^{p}}}$. Then
$$\frac{{{\alpha }_{n}}\left( x \right)}{{{\left( 1-x \right)}^{n+1}}}=\sum\limits_{p=1}^{n}{\frac{{{\alpha }_{p}}{{x}^{p}}}{{{\left( 1-x \right)}^{n+1}}}}=\sum\limits_{m=0}^{\infty }{{{x}^{m}}\sum\limits_{p=1}^{n}{\frac{{{\alpha }_{p}}{}^{\left( p \right)}u_{n}^{{}}\left( m \right)}{n!}}}=\sum\limits_{m=0}^{\infty }{\frac{{{u}_{n}}\left( m \right)}{n!}{{x}^{m}}}.$$
Thus, 
$$U_{2}^{-1}=\left( \begin{matrix}
   1 & -1  \\
   1 & 1  \\
\end{matrix} \right),   \qquad U_{3}^{-1}=\left( \begin{matrix}
   2 & -1 & 2  \\
   3 & 0 & -3  \\
   1 &  1 &   1  \\
\end{matrix} \right),   \qquad U_{4}^{-1}=\left( \begin{matrix}
   6 & -2 &  2 & -6  \\
   11 & -1 & -1 &  11  \\
   6 &  2 & -2 & -6  \\
   1 & 1 & 1 & 1  \\
\end{matrix} \right).$$
{\bfseries Remark 1.} If
$${{\tilde{\alpha }}_{n}}\left( x \right)={{\left( 1-x \right)}^{m}}{{c}_{n-m}}\left( x \right),  \qquad m<n,$$
where ${{c}_{n-m}}\left( x \right)$ is the polynomial of degree $<n-m$, (in this case ${{a}_{1}}=0$ and the matrix ${{\left( 1,\log a\left( x \right) \right)}_{{{e}^{x}}}}$ has no inverse), then, as follows from definition of the transformation $U_{n}^{-1}$, 
$$U_{n}^{-1}{{\left( 1-x \right)}^{m}}{{c}_{n-m}}\left( x \right)=\frac{n!}{\left( n-m \right)!}U_{n-m}^{-1}{{c}_{n-m}}\left( x \right),$$
or
$$U_{n}^{-1}\left( {{\left( 1-x \right)}^{m}},x \right){{I}_{n-m}}=\frac{n!}{\left( n-m \right)!}U_{n-m}^{-1}.$$
For example,
$$\left( \begin{matrix}
   2 & -1 & 2  \\
   3 & 0 & -3  \\
   1 & 1 & 1  \\
\end{matrix} \right)\left( \begin{matrix}
   1 & 0  \\
   -1 & 1  \\
   0 & -1  \\
\end{matrix} \right)=\left( \begin{matrix}
   3 & -3  \\
   3 & 3  \\
\end{matrix} \right)=3U_{2}^{-1}.$$
Accordingly, if ${{c}_{n-m}}\left( x \right)$ is the polynomial of degree $n-m-1$, then
$${{U}_{n}}{{c}_{n-m}}\left( x \right)={{\left( 1-x \right)}^{m}}\frac{\left( n-m \right)!}{n!}{{U}_{n-m}}{{c}_{n-m}}\left( x \right),$$
or
$$\left( {{\left( 1-x \right)}^{-m}},x \right){{U}_{n}}{{I}_{n-m}}=\frac{\left( n-m \right)!}{n!}{{U}_{n-m}}.$$
For example,
$$\left( \begin{matrix}
   1 & 0 & 0  \\
   1 & 1 & 0  \\
   1 & 1 & 1  \\
\end{matrix} \right)\frac{1}{3!}\left( \begin{matrix}
   1 & 1  \\
   -2 & 0  \\
   1 & -1  \\
\end{matrix} \right)=\frac{1}{3!}\left( \begin{matrix}
   1 & 1  \\
   -1 & 1  \\
\end{matrix} \right)=\frac{1}{3}{{U}_{2}}.$$

Denote  ${{\tilde{I}}_{n}}E{{\tilde{I}}_{n}}={{V}_{n}}$. For example,
$${{V}_{2}}=\left( \begin{matrix}
   1 & 0  \\
   1 & 1  \\
\end{matrix} \right),   \qquad{{V}_{3}}=\left( \begin{matrix}
   1 & 0 & 0  \\
   2 & 1 & 0  \\
   1 & 1 & 1  \\
\end{matrix} \right),    \qquad{{V}_{4}}=\left( \begin{matrix}
   1 & 0 & 0 & 0  \\
   3 & 1 & 0 & 0  \\
   3 & 2 & 1 & 0  \\
   1 & 1 & 1 & 1  \\
\end{matrix} \right),$$
$$V_{2}^{-1}=\left( \begin{matrix}
   \text{  }1 & 0  \\
   -1 & 1  \\
\end{matrix} \right),  \qquad V_{3}^{-1}=\left( \begin{matrix}
   1 & 0 & \text{  }0  \\
   -2 & 1 & 0  \\
   1 & -1 & 1  \\
\end{matrix} \right),  \qquad V_{4}^{-1}=\left( \begin{matrix}
  1 & 0 & 0 & 0  \\
   -3 & 1 & 0 & 0  \\
  3 & -2 & 1 & 0  \\
   -1 & 1 & -1 & 1  \\
\end{matrix} \right).$$
If ${{c}_{n}}\left( x \right)$ is the polynomial of degree $<n$, then
$${{V}_{n}}{{c}_{n}}\left( x \right)={{\left( 1+x \right)}^{n-1}}{{c}_{n}}\left( \frac{x}{1+x} \right),  \qquad V_{n}^{-1}{{c}_{n}}\left( x \right)={{\left( 1-x \right)}^{n-1}}{{c}_{n}}\left( \frac{x}{1-x} \right).$$
We have already found out that
$$V_{n}^{-1}{{\tilde{v}}_{n}}\left( x \right)={{\tilde{\alpha }}_{n}}\left( x \right),$$
where
$${{\tilde{v}}_{n}}\left( x \right)=\frac{1}{x}{{v}_{n}}\left( x \right), \qquad{{v}_{n}}\left( x \right)=\left[ n,\to  \right]\left( 1,a\left( x \right)-1 \right),$$
and hence
   $$U_{n}^{-1}V_{n}^{-1}{{\tilde{v}}_{n}}\left( x \right)={{\tilde{u}}_{n}}\left( x \right).$$
It follows from Remark 1, that
$$[\uparrow ,p]\left( U_{n}^{-1}V_{n}^{-1} \right)=U_{n}^{-1}{{\left( 1-x \right)}^{n-p-1}}{{x}^{p}}=\frac{n!}{\left( p+1 \right)!}\frac{1}{x}\prod\limits_{m=0}^{p}{\left( x-m \right)}=$$
$$=\frac{n!}{\left( p+1 \right)!}\sum\limits_{m=1}^{p+1}{s\left( p+1,\text{ }m \right){{x}^{m-1}}},$$
where $s\left( p+1,\text{ }m \right)$ are the Stirling numbers of the first kind. Hence
$$[\uparrow ,p]\left( {{V}_{n}}{{U}_{n}} \right)=\frac{1}{n!}\sum\limits_{m=1}^{p+1}{m!S\left( p+1,\text{ }m \right)}\text{ }{{x}^{m-1}},$$
where $S\left( p+1,\text{ }m \right)$ are the Stirling numbers of the second kind. For example,
$$U_{4}^{-1}V_{4}^{-1}=4!\left( \begin{matrix}
   1 & -1 &2 & -6  \\
   0 & 1 & -3 & 11  \\
   0 & 0 & 1 & -6  \\
   0 & 0 & 0 & 1  \\
\end{matrix} \right)\left( \begin{matrix}
   1 & 0 & 0 & 0  \\
   0 & \frac{1}{2} & 0 & 0  \\
   0 & 0 & \frac{1}{3!} & 0  \\
   0 & 0 & 0 & \frac{1}{4!}  \\
\end{matrix} \right),$$
$${{V}_{4}}{{U}_{4}}=\frac{1}{4!}\left( \begin{matrix}
   1 & 0 & 0 & 0  \\
   0 & 2 & 0 & 0  \\
   0 & 0 & 3! & 0  \\
   0 & 0 & 0 & 4!  \\
\end{matrix} \right)\left( \begin{matrix}
   1 & 1 & 1 & 1  \\
   0 & 1 & 3 & 7  \\
   0 & 0 & 1 & 6  \\
   0 & 0 & 0 & 1  \\
\end{matrix} \right).$$
{\bfseries Remark 2.} Coefficients of the polynomials ${{v}_{n}}\left( x \right)$, ${{u}_{n}}\left( x \right)$  are associated by the relationship:
$${{v}_{n}}\left( x \right)=\sum\limits_{m=1}^{n}{{{B}_{n,m}}\left( {{a}_{1}},{{a}_{2}},...,{{a}_{n}} \right){{x}^{m}}},  \qquad{{u}_{n}}\left( x \right)=n!\sum\limits_{m=1}^{n}{\frac{{{B}_{n,m}}\left( {{b}_{1}},{{b}_{2}},...,{{b}_{n}} \right)}{m!}{{x}^{m}}},$$
where
$${{b}_{p}}=\left[ {{x}^{p}} \right]\log a\left( x \right)=\sum\limits_{m=1}^{p}{{{\left( -1 \right)}^{m+1}}\frac{{{B}_{p,m}}\left( {{a}_{1}},{{a}_{2}},...,{{a}_{p}} \right)}{m}},$$
$${{B}_{n,m}}\left( {{a}_{1}},{{a}_{2}},...,{{a}_{n}} \right)=\sum{\frac{m!}{{{m}_{1}}!{{m}_{2}}!\text{ }...\text{ }{{m}_{n}}!}}\text{ }a_{1}^{{{m}_{1}}}a_{2}^{{{m}_{2}}}...\text{ }a_{n}^{{{m}_{n}}},$$   
$${{B}_{n,m}}\left( {{b}_{1}},{{b}_{2}},...,{{b}_{n}} \right)=\sum{\frac{m!}{{{m}_{1}}!{{m}_{2}}!\text{ }...\text{ }{{m}_{n}}!}\text{ }}b_{1}^{{{m}_{1}}}b_{2}^{{{m}_{2}}}...\text{ }b_{n}^{{{m}_{n}}},$$
expressions $\prod\nolimits_{p=1}^{n}{a_{p}^{{{m}_{p}}}}$, $\prod\nolimits_{p=1}^{n}{b_{p}^{{{m}_{p}}}}$ corresponding to the partition $n=\sum\nolimits_{p=1}^{n}{p{{m}_{p}}}$, $\sum\nolimits_{p=1}^{n}{{{m}_{p}}}=m$, and summation is done over all partitions of number $n$ to $m$ parts.
\section{ Examples}
{\bfseries Example 1.}
$$a\left( x \right)=\frac{1+x}{1-x},   \qquad a\left( x \right)-1=\frac{2x}{1-x},  \qquad{{\tilde{v}}_{n}}\left( x \right)={{2}^{n}}{{\left( \frac{1}{2}+x \right)}^{n-1}},$$
$${{\tilde{\alpha }}_{n}}\left( x \right)=V_{n}^{-1}{{2}^{n}}{{\left( \frac{1}{2}+x \right)}^{n-1}}=2{{\left( 1+x \right)}^{n-1}},$$
$${{u}_{n}}\left( x \right)=xU_{n}^{-1}{{\tilde{\alpha }}_{n}}\left( x \right)=2\sum\limits_{p=0}^{n-1}{\left( \begin{matrix}
   n-1  \\
   p  \\
\end{matrix} \right)\prod\limits_{m=0}^{n-1}{\left( x-p+m \right)}}=$$
$$=xU_{n}^{-1}V_{n}^{-1}{{\tilde{v}}_{n}}\left( x \right)=n!\sum\limits_{p=0}^{n-1}{\left( \begin{matrix}
   n-1  \\
   p  \\
\end{matrix} \right)\frac{{{2}^{p+1}}}{\left( p+1 \right)!}}\prod\limits_{m=0}^{p}{\left( x-m \right)}.$$
{\bfseries Example 2.} Let
$$b\left( xg\left( x \right) \right)=g\left( x \right),    \qquad g\left( x{{b}^{-1}}\left( x \right) \right)=b\left( x \right).$$
Then by the Lagrange inversion theorem 
$$\left[ {{x}^{n}} \right]{{g}^{m}}\left( x \right)=\left[ {{x}^{n}} \right]\left( 1-x{{\left( \log b\left( x \right) \right)}^{\prime }} \right){{b}^{m+n}}\left( x \right),$$
$$\left( 1,xg\left( x \right) \right)=\left( \begin{matrix}
   g_{0}^{0} & 0 & 0 & 0 & \cdots   \\
   g_{1}^{0} & g_{0}^{1} & 0 & 0 & \cdots   \\
   g_{2}^{0} & g_{1}^{1} & g_{0}^{2} & 0 & \cdots   \\
   g_{3}^{0} & g_{2}^{1} & g_{1}^{2} & g_{0}^{3} & \cdots   \\
   \vdots  & \vdots  & \vdots  & \vdots  & \ddots   \\
\end{matrix} \right)=\left( \begin{matrix}
   c_{0}^{0} & 0 & 0 & 0 & \cdots   \\
   c_{1}^{1} & c_{0}^{1} & 0 & 0 & \cdots   \\
   c_{2}^{2} & c_{1}^{2} & c_{0}^{2} & 0 & \cdots   \\
   c_{3}^{3} & c_{2}^{3} & c_{1}^{3} & c_{0}^{3} & \cdots   \\
   \vdots  & \vdots  & \vdots  & \vdots  & \ddots   \\
\end{matrix} \right),$$
where
  $$g_{n}^{m}=\left[ {{x}^{n}} \right]{{g}^{m}}\left( x \right),  \qquad c_{n}^{m}=\left[ {{x}^{n}} \right]\left( 1-x{{\left( \log b\left( x \right) \right)}^{\prime }} \right){{b}^{m}}\left( x \right),$$
i.e.
$$\left[ n,\to  \right]\left( 1,xg\left( x \right) \right)=\left[ n,\to  \right]\left( \left( 1-x{{\left( \log b\left( x \right) \right)}^{\prime }} \right){{b}^{n}}\left( x \right),x \right).$$
If
$$g\left( x \right)=\frac{x}{2}+{{\left( 1+\frac{{{x}^{2}}}{4} \right)}^{{1}/{2}\;}},$$
then
$$b\left( x \right)={{\left( 1+x \right)}^{{1}/{2}\;}},    \qquad1-x{{\left( \log b\left( x \right) \right)}^{\prime }}=\left( 1+\frac{x}{2} \right){{\left( 1+x \right)}^{-1}},$$
$$\left[ 2n,\to  \right]\left( 1,xg\left( x \right) \right)=\left( \frac{1}{2}+x \right){{x}^{n}}{{\left( 1+x \right)}^{n-1}},  \qquad n>0,$$
Let
$$a\left( x \right)={{\left( \frac{x}{2}+{{\left( 1+\frac{{{x}^{2}}}{4} \right)}^{{1}/{2}\;}} \right)}^{2}}.$$
Then
$$a\left( x \right)-1=x{{a}^{{1}/{2}\;}}\left( x \right)=x\left( \frac{x}{2}+{{\left( 1+\frac{{{x}^{2}}}{4} \right)}^{{1}/{2}\;}} \right),$$
$${{\tilde{v}}_{2n}}\left( x \right)=\left( \frac{1}{2}+x \right){{x}^{n-1}}{{\left( 1+x \right)}^{n-1}},$$
$${{\tilde{\alpha }}_{2n}}\left( x \right)=V_{2n}^{-1}\left( \frac{1}{2}+x \right){{x}^{n-1}}{{\left( 1+x \right)}^{n-1}}=\frac{1}{2}\left( 1+x \right){{x}^{n-1}};$$
$${{u}_{2n}}\left( x \right)=xU_{2n}^{-1}{{\tilde{\alpha }}_{2n}}\left( x \right)=\frac{1}{2}\prod\limits_{m=0}^{2n-1}{\left( x-n+1+m \right)+\frac{1}{2}\prod\limits_{m=0}^{2n-1}{\left( x-n+m \right)}}=$$
$$=x\prod\limits_{m=1}^{2n-1}{\left( x+n-m \right)}=\prod\limits_{m=0}^{n-1}{\left( {{x}^{2}}-{{m}^{2}} \right)}.$$
{\bfseries Example 3.} This example was considered in [16]. We will replace the convolution arrays by the Riordan arrays. Denote
$$\left[ n,\to  \right]\left( \frac{1}{1-x-k{{x}^{2}}},\frac{1}{1-x-k{{x}^{2}}} \right)=\frac{{{N}_{n}}\left( x \right)}{{{\left( 1-x \right)}^{n+1}}},$$
$$\left[ n,\to  \right]\left( \frac{1}{1-kx-{{x}^{2}}},\frac{1}{1-kx-{{x}^{2}}} \right)=\frac{N_{n}^{*}\left( x \right)}{{{\left( 1-x \right)}^{n+1}}}.$$
Then
$${{N}_{n}}\left( x \right)=\left[ n,\to  \right]\left( \frac{1}{1-x-k{{x}^{2}}},\frac{-k{{x}^{2}}}{1-x-k{{x}^{2}}} \right),$$
$$N_{n}^{*}\left( x \right)=\left[ n,\to  \right]\left( \frac{1}{1-kx-{{x}^{2}}},\frac{-{{x}^{2}}}{1-kx-{{x}^{2}}} \right).$$
For example,
$$\left( \frac{1}{1-x-{{x}^{2}}},\frac{1}{1-x-{{x}^{2}}} \right)=\left( \begin{matrix}
   1 & 1 & 1 & 1 & 1 & \cdots   \\
   1 & 2 & 3 & 4 & 5 & \cdots   \\
   2 & 5 & 9 & 14 & 20 & \cdots   \\
   3 & 10 & 22 & 40 & 65 & \cdots   \\
   5 & 20 & 51 & 105 & 190 & \cdots   \\
   \vdots  & \vdots  & \vdots  & \vdots  & \vdots  & \ddots   \\
\end{matrix} \right),$$
$$\left( \frac{1}{1-x-{{x}^{2}}},\frac{-{{x}^{2}}}{1-x-{{x}^{2}}} \right)=\left( \begin{matrix}
   1 & 0 & 0 & 0 & \cdots   \\
   1 & 0 & 0 & 0 & \cdots   \\
   2 & -1 & 0 & 0 & \cdots   \\
   3 & -2 & 0 & 0 & \cdots   \\
   5 & -5 & 1 & 0 & \cdots   \\
   8 & -10 & 3 & 0 & \cdots   \\
   \vdots  & \vdots  & \vdots  & \vdots  & \ddots   \\
\end{matrix} \right),$$
$$\frac{2-x}{{{\left( 1-x \right)}^{3}}}=2+5x+9{{x}^{2}}+14{{x}^{3}}+... ,$$ 
$$\frac{3-2x}{{{\left( 1-x \right)}^{4}}}=3+10x+22{{x}^{2}}+40{{x}^{3}}+... ,$$
$$\frac{5-5x+{{x}^{2}}}{{{\left( 1-x \right)}^{5}}}=5+20x+51{{x}^{2}}+105{{x}^{3}}+... .$$
We generalize this example using the transformation $V_{n}^{-1}$. Let $a\left( x \right)={{\left( 1+\varphi x+\beta {{x}^{2}} \right)}^{-1}}$. Then
$$\left( 1,{{a}^{-1}}\left( x \right)-1 \right)=\left( \begin{matrix}
   1 & \text{ }0 & 0 & 0 & 0 & \cdots   \\
   0 & \text{ }\varphi  & 0 & 0 & 0 & \cdots   \\
   0 & \text{ }\beta  & {{\varphi }^{2}} & 0 & 0 & \cdots   \\
   0 & \text{ }0 & 2\varphi \beta  & {{\varphi }^{3}} & 0 & \cdots   \\
   0 & \text{ }0 & {{\beta }^{2}} & 3{{\varphi }^{2}}\beta  & \text{ }{{\varphi }^{4}} & \cdots   \\
   \vdots  & \text{ }\vdots  & \vdots  & \vdots  & \vdots  & \ddots   \\
\end{matrix} \right).$$
Applying the transformation ${{\hat{I}}_{n}}$ to the $n$th row of this matrix, we obtain the matrix
$$\left( \frac{1}{1-\varphi x},\frac{\beta {{x}^{2}}}{1-\varphi x} \right)=\left( \begin{matrix}
   1 & 0 & 0 & \cdots   \\
   \varphi  & 0 & 0 & \cdots   \\
   {{\varphi }^{2}} & \beta  & 0 & \cdots   \\
   {{\varphi }^{3}} & 2\varphi \beta  & 0 & \cdots   \\
   {{\varphi }^{4}} & 3{{\varphi }^{2}}\beta  & \text{ }{{\beta }^{2}} & \cdots   \\
   \vdots  & \vdots  & \vdots  & \ddots   \\
\end{matrix} \right).$$
Since
$${{\tilde{\alpha }}_{n}}\left( x \right)={{\left( -1 \right)}^{n}}{{\tilde{I}}_{n}}V_{n}^{-1}\tilde{v}_{n}^{\left( -1 \right)}\left( x \right)={{\left( -1 \right)}^{n}}{{E}^{-1}}{{\tilde{I}}_{n}}\tilde{v}_{n}^{\left( -1 \right)}\left( x \right),$$ 
$${{\tilde{I}}_{n}}\tilde{v}_{n}^{\left( -1 \right)}\left( x \right)={{\hat{I}}_{n}}v_{n}^{\left( -1 \right)}\left( x \right),$$
then the polynomial ${{\tilde{\alpha }}_{n}}\left( x \right)$ corresponds to the $n$th row of the matrix
$$\left( 1,-x \right)\left( \frac{1}{1-\varphi x},\frac{\beta {{x}^{2}}}{1-\varphi x} \right)\left( \frac{1}{1+x},\frac{x}{1+x} \right)=\left( \frac{1}{1+\varphi x+\beta {{x}^{2}}},\frac{\beta {{x}^{2}}}{1+\varphi x+\beta {{x}^{2}}} \right).$$
Really,
$$\sum\limits_{n=1}^{\infty }{{{{\tilde{\alpha }}}_{n}}}\left( t \right){{x}^{n}}=\frac{-\varphi x-\beta \left( 1-t \right){{x}^{2}}}{1+\varphi x+\beta \left( 1-t \right){{x}^{2}}},$$
$$\sum\limits_{n=0}^{\infty }{{{\alpha }_{n}}}\left( t \right){{x}^{n}}=1+t\sum\limits_{n=1}^{\infty }{{{{\tilde{\alpha }}}_{n}}}\left( t \right){{x}^{n}}=\frac{1+\varphi \left( 1-t \right)x+\beta {{\left( 1-t \right)}^{2}}{{x}^{2}}}{1+\varphi x+\beta \left( 1-t \right){{x}^{2}}},$$ 
which corresponds to the formula (2). Since
$$\left[ n,\to  \right]\left( \frac{1}{1+\varphi x},\frac{\beta {{x}^{2}}}{1+\varphi x} \right)={{r}_{n}}\prod\limits_{m=1}^{\left\lfloor {n}/{2}\; \right\rfloor }{\left( \beta x+\frac{{{\varphi }^{2}}}{4}{{\sec }^{2}}\frac{m}{n+1}\pi  \right)},$$
where ${{r}_{2p}}=1$, ${{r}_{2p+1}}=-\left( p+1 \right)\varphi $, (i.e. these polynomials are associated in a certain way with the Chebyshev polynomials), then
$${{\tilde{\alpha }}_{n}}\left( x \right)={{r}_{n}}\prod\limits_{m=1}^{\left\lfloor {n}/{2}\; \right\rfloor }{\left( \beta x+\frac{{{\left( {\varphi }/{2}\; \right)}^{2}}-\beta {{\cos }^{2}}\frac{m}{n+1}\pi }{{{\cos }^{2}}\frac{m}{n+1}\pi } \right)}.$$
{\bfseries Example 4.}  In [24] Carlitz and Hoggatt considered the following generalization of Euler polynomials:
$$G_{n}^{\left( p \right)}\left( x \right)={{\left( 1-x \right)}^{pn+1}}\sum\limits_{m=0}^{\infty }{{{\left( \begin{matrix}
   m+p-1  \\
   p  \\
\end{matrix} \right)}^{n}}}{{x}^{m}},$$
$G_{n}^{\left( 1 \right)}\left( x \right)={{A}_{n}}\left( x \right)$, $G_{n}^{\left( p \right)}\left( x \right)$ is the polynomial of degree $pn-p+1$, such that
$$\left[ {{x}^{m}} \right]G_{n}^{\left( p \right)}\left( x \right)=\left[ {{x}^{pn-p-m+2}} \right]G_{n}^{\left( p \right)}\left( x \right),  \qquad1\le m\le pn-p+1;$$
$$G_{n}^{\left( p \right)}\left( 1 \right)=\frac{\left( pn \right)!}{{{\left( p! \right)}^{n}}}.$$
Properties of these polynomials will become more transparent if we associate them with the following construction. We will consider the formal Dirichlet series $a\left( s \right)=\sum\nolimits_{n=1}^{\infty }{{{{a}_{n}}}/{{{n}^{s}}}\;}$ as the generating function of the sequence ${{\left( {{a}_{n}} \right)}_{n\ge 0}}$, ${{a}_{0}}=0$. 
Matrix whose $n$th column has the generating function ${{a}^{n}}\left( s \right)$ is denoted $\left\langle a\left( s \right) \right\rangle $.
For example, for the Riemann zeta function: 
$$\left\langle \zeta \left( s \right) \right\rangle =\left( \begin{matrix}
   0 & 0 & 0 & 0 & 0 & \cdots   \\
   1 & 1 & 1 & 1 & 1 & \cdots   \\
   0 & 1 & 2 & 3 & 4 & \cdots   \\
   0 & 1 & 2 & 3 & 4 & \cdots   \\
   0 & 1 & 3 & 6 & 10 & \cdots   \\
   0 & 1 & 2 & 3 & 4 & \cdots   \\
   0 & 1 & 4 & 9 & 16 & \cdots   \\
   0 & 1 & 2 & 3 & 4 & \cdots   \\
   0 & 1 & 4 & 10 & 20 & \cdots   \\
   0 & 1 & 3 & 6 & 10 & \cdots   \\
   0 & 1 & 4 & 9 & 16 & \cdots   \\
   0 & 1 & 2 & 3 & 4 & \cdots   \\
   0 & 1 & 6 & 18 & 40 & \cdots   \\
   \vdots  & \vdots  & \vdots  & \vdots  & \vdots  & \ddots   \\
\end{matrix} \right),  \qquad\left\langle {{\zeta }^{-1}}\left( s \right) \right\rangle =\left( \begin{matrix}
   0 & 0 & 0 & 0 & 0 & \cdots   \\
   1 & 1 & 1 & 1 & 1 & \cdots   \\
   0 & -1 & -2 & -3 & -4 & \cdots   \\
   0 & -1 & -2 & -3 & -4 & \cdots   \\
   0 & 0 & 1 & 3 & 6 & \cdots   \\
   0 & -1 & -2 & -3 & -4 & \cdots   \\
   0 & 1 & 4 & 9 & 16 & \cdots   \\
   0 & -1 & -2 & -3 & -4 & \cdots   \\
   0 & 0 & 0 & -1 & -4 & \cdots   \\
   0 & 0 & 1 & 3 & 6 & \cdots   \\
   0 & 1 & 4 & 9 & 16 & \cdots   \\
   0 & -1 & -2 & -3 & -4 & \cdots   \\
   0 & 0 & -2 & -9 & -24 & \cdots   \\
   \vdots  & \vdots  & \vdots  & \vdots  & \vdots  & \ddots   \\
\end{matrix} \right),$$
$$\left\langle \zeta \left( s \right)-1 \right\rangle =\left( \begin{matrix}
   0 & 0 & 0 & 0 & \cdots   \\
   1 & 0 & 0 & 0 & \cdots   \\
   0 & 1 & 0 & 0 & \cdots   \\
   0 & 1 & 0 & 0 & \cdots   \\
   0 & 1 & 1 & 0 & \cdots   \\
   0 & 1 & 0 & 0 & \cdots   \\
   0 & 1 & 2 & 0 & \cdots   \\
   0 & 1 & 0 & 0 & \cdots   \\
   0 & 1 & 2 & 1 & \cdots   \\
   0 & 1 & 1 & 0 & \cdots   \\
   0 & 1 & 2 & 0 & \cdots   \\
   0 & 1 & 0 & 0 & \cdots   \\
   0 & 1 & 4 & 3 & \cdots   \\
   \vdots  & \vdots  & \vdots  & \vdots  & \ddots   \\
\end{matrix} \right),   \qquad\left\langle \log \zeta \left( s \right) \right\rangle =\left( \begin{matrix}
   0 & 0 & 0 & 0 & \cdots   \\
   1 & 0 & 0 & 0 & \cdots   \\
   0 & 1 & 0 & 0 & \cdots   \\
   0 & 1 & 0 & 0 & \cdots   \\
   0 & {1}/{2}\; & 1 & 0 & \cdots   \\
   0 & 1 & 0 & 0 & \cdots   \\
   0 & 0 & 2 & 0 & \cdots   \\
   0 & 1 & 0 & 0 & \cdots   \\
   0 & {1}/{3}\; & 1 & 1 & \cdots   \\
   0 & {1}/{2}\; & 1 & 0 & \cdots   \\
   0 & 0 & 2 & 0 & \cdots   \\
   0 & 1 & 0 & 0 & \cdots   \\
   0 & 0 & 1 & 3 & \cdots   \\
   \vdots  & \vdots  & \vdots  & \vdots  & \ddots   \\
\end{matrix} \right).$$
Such matrices are considered in [25]. If the matrix $\left\langle a\left( s \right) \right\rangle $, ${{a}_{1}}=0$, is multiplied from the right by the Riordan matrix $\left( 1,g\left( x \right) \right)$, $g\left( x \right)=\sum\nolimits_{n=0}^{\infty }{{{g}_{n}}{{x}^{n}}}$,  the result is the matrix $\left\langle \sum\nolimits_{n=0}^{\infty }{{{g}_{n}}{{a}^{n}}\left( s \right)} \right\rangle $. In particular,  if ${{a}_{1}}=1$,
$$\left\langle a\left( s \right)-1 \right\rangle \left( 1,1+x \right)=\left\langle a\left( s \right) \right\rangle ,  \qquad\left\langle a\left( s \right)-1 \right\rangle \left( 1,\frac{1}{1+x} \right)=\left\langle {{a}^{-1}}\left( s \right) \right\rangle ,$$
$$\left\langle \log a\left( s \right) \right\rangle \left( 1,{{e}^{x}} \right)=\left\langle a\left( s \right) \right\rangle .$$
We associate rows of the matrices $\left\langle a\left( s \right) \right\rangle $ with the formal power series, which are the generating functions of their elements. For polynomials similar to polynomials associated with the GEP, we use the same notation. Then ($n>1)$
$$\left[ n,\to  \right]\left\langle a\left( s \right)-1 \right\rangle ={{v}_{n}}\left( x \right),  \qquad\frac{1}{x}{{v}_{n}}\left( x \right)={{\tilde{v}}_{n}}\left( x \right),$$
$$\left[ n,\to  \right]\left\langle a\left( s \right) \right\rangle =\frac{{{\alpha }_{n}}\left( x \right)}{{{\left( 1-x \right)}^{v\left( n \right)+1}}} ,  \qquad\left[ n,\to  \right]\left\langle {{a}^{-1}}\left( s \right) \right\rangle =\frac{\alpha _{n}^{\left( -1 \right)}\left( x \right)}{{{\left( 1-x \right)}^{v\left( n \right)+1}}},$$
where
$$\alpha _{n}^{\left( -1 \right)}\left( x \right)={{\left( -1 \right)}^{v\left( n \right)}}x{{\hat{I}}_{v\left( n \right)}}{{\alpha }_{n}}\left( x \right),  \qquad{{\alpha }_{n}}\left( x \right)=xV_{v\left( n \right)}^{-1}{{\tilde{v}}_{n}}\left( x \right),$$
$v\left( n \right)$ is the degree of polynomial ${{v}_{n}}\left( x \right)$;
$$\left[ n,\to  \right]\left( {{\left| {{e}^{x}} \right|}^{-1}}\left\langle \log a\left( s \right) \right\rangle \left| {{e}^{x}} \right| \right)={{u}_{n}}\left( x \right),  \qquad\frac{1}{x}{{u}_{n}}\left( x \right)={{\tilde{u}}_{n}}\left( x \right),$$
$${{a}^{m}}\left( s \right)=\sum\limits_{n=0}^{\infty }{\frac{{{u}_{n}}\left( m \right)}{n!{{n}^{s}}}},   \qquad\frac{{{\alpha }_{n}}\left( x \right)}{{{\left( 1-x \right)}^{u\left( n \right)+1}}}=\sum\limits_{m=0}^{\infty }{\frac{{{u}_{n}}\left( m \right)}{n!}}{{x}^{m}},$$
$${{\alpha }_{n}}\left( x \right)=x\frac{u\left( n \right)!}{n!}{{U}_{u\left( n \right)}}{{\tilde{u}}_{n}}\left( x \right),$$
where $u\left( n \right)$ is the degree of polynomial ${{u}_{n}}\left( x \right)$, equal to the degree of polynomial ${{v}_{n}}\left( x \right)$. Matrix $\left\langle a\left( s \right)-1 \right\rangle $ has the form:
$$\left\langle a\left( s \right)-1 \right\rangle =\left( \begin{matrix}
   0 & 0 & 0 & 0 & 0 & \ldots   \\
   1 & 0 & 0 & 0 & 0 & \ldots   \\
   0 & {{a}_{2}} & 0 & 0 & 0 & \ldots   \\
   0 & {{a}_{3}} & 0 & 0 & 0 & \ldots   \\
   0 & {{a}_{4}} & a_{2}^{2} & 0 & 0 & \ldots   \\
   0 & {{a}_{5}} & 0 & 0 & 0 & \ldots   \\
   0 & {{a}_{6}} & 2{{a}_{2}}{{a}_{3}} & 0 & 0 & \ldots   \\
   0 & {{a}_{7}} & 0 & 0 & 0 & \ldots   \\
   0 & {{a}_{8}} & 2{{a}_{2}}{{a}_{4}} & a_{2}^{3} & 0 & \ldots   \\
   0 & {{a}_{9}} & a_{3}^{2} & 0 & 0 & \ldots   \\
   0 & {{a}_{10}} & 2{{a}_{2}}{{a}_{5}} & 0 & 0 & \ldots   \\
   0 & {{a}_{11}} & 0 & 0 & 0 & \ldots   \\
   0 & {{a}_{12}} & 2{{a}_{2}}{{a}_{6}}+2{{a}_{4}}{{a}_{3}} & 3a_{2}^{2}{{a}_{3}} & 0 & \ldots   \\
   0 & {{a}_{13}} & 0 & 0 & 0 & \ldots   \\
   0 & {{a}_{14}} & 2{{a}_{2}}{{a}_{7}} & 0 & 0 & \ldots   \\
   0 & {{a}_{15}} & 2{{a}_{3}}{{a}_{5}} & 0 & 0 & \ldots   \\
   0 & {{a}_{16}} & 2{{a}_{2}}{{a}_{8}}+a_{4}^{2} & 3a_{2}^{2}{{a}_{4}} & a_{2}^{4} & \ldots   \\
   \vdots  & \vdots  & \vdots  & \vdots  & \vdots  & \ddots   \\
\end{matrix} \right)$$
$${{v}_{n}}\left( x \right)=\sum\limits_{m=1}^{v\left( n \right)}{{{{\tilde{B}}}_{n,m}}\left( {{a}_{2}},{{a}_{3}},...,{{a}_{n}} \right){{x}^{m}}},  \qquad n>1,$$
where
$${{\tilde{B}}_{n,m}}\left( {{a}_{2}},{{a}_{3}},...,{{a}_{n}} \right)=\sum{\frac{m!}{{{m}_{2}}!{{m}_{3}}!\text{ }...\text{ }{{m}_{n}}!}}\text{ }a_{2}^{{{m}_{2}}}a_{3}^{{{m}_{3}}}...\text{ }a_{n}^{{{m}_{n}}},$$  
expression $\prod\nolimits_{p=2}^{n}{a_{p}^{{{m}_{p}}}}$ corresponding to the decomposition $n=\prod\nolimits_{p=2}^{n}{{{p}^{{{m}_{p}}}}}$, $\sum\nolimits_{p=2}^{n}{{{m}_{p}}}=m$, and summation is done over all decompositions of number $n$ into $m$ factors . If $\log a\left( s \right)=b\left( s \right)$, then
$${{b}_{p}}=\sum\limits_{m=1}^{v\left( p \right)}{{{\left( -1 \right)}^{m+1}}\frac{{{{\tilde{B}}}_{p,m}}\left( {{a}_{2}},{{a}_{3}},...,{{a}_{p}} \right)}{m}},  \qquad{{u}_{n}}\left( x \right)=n!\sum\limits_{m=1}^{u\left( n \right)}{\frac{{{{\tilde{B}}}_{n,m}}\left( {{b}_{2}},{{b}_{3}},...,{{b}_{n}} \right)}{m!}{{x}^{m}}}.$$

If $a\left( s \right)=\zeta \left( s \right)$, then
$${{u}_{0}}\left( x \right)=0,   \qquad{{u}_{1}}\left( x \right)=1,    \qquad\frac{{{u}_{n}}\left( x \right)}{n!}=\frac{{{\left( x \right)}^{{{m}_{1}}}}{{\left( x \right)}^{{{m}_{2}}}}...{{\left( x \right)}^{{{m}_{r}}}}}{{{m}_{1}}!{{m}_{2}}!\text{ }...\text{ }{{m}_{r}}!},$$
where
$${{\left( x \right)}^{{{m}_{i}}}}=x\left( x+1 \right)\left( x+2 \right)...\left( x+{{m}_{i}}-1 \right),$$
 $n=p_{1}^{{{m}_{1}}}p_{2}^{{{m}_{2}}}...\text{ }p_{r}^{{{m}_{r}}}$ is the canonical decomposition of number $n$. If ${{m}_{1}}={{m}_{2}}=...={{m}_{r}}=p$, then
$$\frac{{{u}_{n}}\left( x \right)}{n!}={{\left( \frac{{{\left( x \right)}^{p}}}{p!} \right)}^{r}}, \qquad u\left( n \right)=pr,  \qquad\left[ {{x}^{pr}} \right]{{u}_{n}}\left( x \right)=\frac{n!}{{{\left( p! \right)}^{r}}},$$
$${{\alpha }_{n}}\left( x \right)=G_{r}^{\left( p \right)}\left( x \right),  \qquad\alpha _{n}^{\left( -1 \right)}\left( x \right)={{\left( -1 \right)}^{pr}}{{x}^{p-1}}G_{r}^{\left( p \right)}\left( x \right).$$
It is clear from this that the sum of coefficients and the degree of the polynomial $G_{r}^{\left( p \right)}\left( x \right)$ can be  defined from the transformations
$$G_{r}^{\left( p \right)}\left( x \right)=x\frac{\left( pr \right)!}{n!}{{U}_{pr}}{{\tilde{u}}_{n}}\left( x \right),  \qquad{{x}^{p-1}}G_{r}^{\left( p \right)}\left( x \right)=x{{\hat{I}}_{pr}}G_{r}^{\left( p \right)}\left( x \right).$$
\section{ GEP and multinomial coefficients}

Denote
$$\left[ n,\to  \right]\left( 1,{{a}^{m}}\left( x \right) \right)=\frac{\alpha _{n}^{\left( m \right)}\left( x \right)}{{{\left( 1-x \right)}^{n+1}}},  \qquad\alpha _{n}^{\left( 1 \right)}\left( x \right)=\alpha _{n}^{{}}\left( x \right),  \qquad\frac{1}{x}\alpha _{n}^{\left( m \right)}\left( x \right)=\tilde{\alpha }_{n}^{\left( m \right)}\left( x \right).$$
Then
$${{U}_{n}}m{{\tilde{u}}_{n}}\left( mx \right)=\tilde{\alpha }_{n}^{\left( m \right)}\left( x \right),\qquad
{{U}_{n}}\left( m,mx \right)U_{n}^{-1}{{\tilde{\alpha }}_{n}}\left( x \right)=\tilde{\alpha }_{n}^{\left( m \right)}\left( x \right).$$
Denote
$${{W}_{\left( n,\text{ }m \right)}}={{U}_{n}}\left( m,mx \right)U_{n}^{-1}.$$
Construct the matrix ${{\left( b\left( x \right),x \right)}_{m}}$ by the rule
$$\left[ n,\to  \right]{{\left( b\left( x \right),x \right)}_{m}}=\left[ mn+m-1,\to  \right]\left( b\left( x \right),x \right).$$
For example,
$${{\left( b\left( x \right),x \right)}_{2}}=\left( \begin{matrix}
   {{b}_{1}} & {{b}_{0}} & 0 & 0 & \cdots   \\
   {{b}_{3}} & {{b}_{2}} & {{b}_{1}} & {{b}_{0}} & \cdots   \\
   {{b}_{5}} & {{b}_{4}} & {{b}_{3}} & {{b}_{2}} & \cdots   \\
   {{b}_{7}} & {{b}_{6}} & {{b}_{5}} & {{b}_{4}} & \cdots   \\
   \vdots  & \vdots  & \vdots  & \vdots  & \ddots   \\
\end{matrix} \right),  \qquad{{\left( b\left( x \right),x \right)}_{3}}=\left( \begin{matrix}
   {{b}_{2}} & {{b}_{1}} & {{b}_{0}} & 0 & \cdots   \\
   {{b}_{5}} & {{b}_{4}} & {{b}_{3}} & {{b}_{2}} & \cdots   \\
   {{b}_{8}} & {{b}_{7}} & {{b}_{6}} & {{b}_{5}} & \cdots   \\
   {{b}_{11}} & {{b}_{10}} & {{b}_{9}} & {{b}_{8}} & \cdots   \\
   \vdots  & \vdots  & \vdots  & \vdots  & \ddots   \\
\end{matrix} \right).$$
{\bfseries Theorem 4.}
$${{W}_{\left( n,\text{ }m \right)}}={{\left( w_{m}^{n+1}\left( x \right),x \right)}_{m}}{{I}_{n}} ,  \qquad w_{m}^{n+1}\left( x \right)={{\left( \frac{1-{{x}^{m}}}{1-x} \right)}^{n+1}}.$$
{\bfseries Proof.} Since
$$\frac{{{{\tilde{\alpha }}}_{n}}\left( x \right)}{{{\left( 1-x \right)}^{n+1}}}=\frac{w_{m}^{n+1}\left( x \right){{{\tilde{\alpha }}}_{n}}\left( x \right)}{{{\left( 1-{{x}^{m}} \right)}^{n+1}}}=\sum\limits_{r=0}^{m-1}{\frac{{{x}^{r}}{{c}_{r}}\left( x \right)}{{{\left( 1-{{x}^{m}} \right)}^{n+1}}}},$$
where
$${{c}_{r}}\left( x \right)=\sum\limits_{p=0}^{\infty }{\left( \left[ {{x}^{mp+r}} \right]w_{m}^{n+1}\left( x \right){{{\tilde{\alpha }}}_{n}}\left( x \right) \right){{x}^{mp}}},$$
and since
$$\left[ {{x}^{p}} \right]\frac{\tilde{\alpha }_{n}^{\left( m \right)}\left( x \right)}{{{\left( 1-x \right)}^{n+1}}}=\left[ {{x}^{mp+m-1}} \right]\frac{{{{\tilde{\alpha }}}_{n}}\left( x \right)}{{{\left( 1-x \right)}^{n+1}}},$$
then
$$\frac{\tilde{\alpha }_{n}^{\left( m \right)}\left( {{x}^{m}} \right)}{{{\left( 1-{{x}^{m}} \right)}^{n+1}}}=\frac{{{c}_{m-1}}\left( x \right)}{{{\left( 1-{{x}^{m}} \right)}^{n+1}}},$$
$$\left[ {{x}^{p}} \right]\tilde{\alpha }_{n}^{\left( m \right)}\left( x \right)=\left[ {{x}^{mp+m-1}} \right]w_{m}^{n+1}\left( x \right){{\tilde{\alpha }}_{n}}\left( x \right),$$
or
$$\tilde{\alpha }_{n}^{\left( m \right)}\left( x \right)={{\left( w_{m}^{n+1}\left( x \right),x \right)}_{m}}{{\tilde{\alpha }}_{n}}\left( x \right).$$
For example (${{\left( w_{m}^{n}\left( x \right) \right)}_{i}}$ means the sequence of coefficients of the polynomial $w_{m}^{n}\left( x \right)$):
$${{\left( w_{2}^{2}\left( x \right) \right)}_{i}}=\left( 1,\text{ }2,\text{ }1 \right), \quad{{\left( w_{2}^{3}\left( x \right) \right)}_{i}}=\left( 1,\text{ }3,\text{ }3,\text{ }1 \right), \quad{{\left( w_{2}^{3}\left( x \right) \right)}_{i}}=\left( 1,\text{ }4,\text{ }6,\text{ }4,\text{ }1 \right);$$
$${{W}_{(1,\text{ }2)}}=\left( 2 \right), \quad{{W}_{(2,\text{ }2)}}=\left( \begin{matrix}
   3 & 1  \\
   1 & 3  \\
\end{matrix} \right),\quad{{W}_{(3,\text{ }2)}}=\left( \begin{matrix}
   4 & 1 & 0  \\
   4 & 6 & 4  \\
   0 & 1 & 4  \\
\end{matrix} \right).$$
$${{\left( w_{3}^{2}\left( x \right) \right)}_{i}}=\left( 1,\text{ }2,\text{ }3,\text{ }2,\text{ }1 \right),   \qquad{{\left( w_{3}^{3}\left( x \right) \right)}_{i}}=\left( 1,\text{ }3,\text{ }6,\text{ }7,\text{ }6,\text{ }3,\text{ }1 \right),$$
$${{\left( w_{3}^{4}\left( x \right) \right)}_{i}}=\left( 1,\text{ }4,\text{ }10,\text{ }16,\text{ }19,\text{ }16,\text{ }10,\text{ }4,\text{ }1 \right);$$
$${{W}_{(1,\text{ }3)}}=\left( 3 \right),\qquad {{W}_{(2,\text{ }3)}}=\left( \begin{matrix}
   6 & 3  \\
   3 & 6  \\
\end{matrix} \right) , \qquad{{W}_{(3,\text{ }3)}}=\left( \begin{matrix}
   10 & 4 & 1  \\
   16 & 19 & 16  \\
   1 & 4 & 10  \\
\end{matrix} \right).$$
$${{\left( w_{4}^{2}\left( x \right) \right)}_{i}}=\left( 1,\text{ }2,\text{ }3,\text{ }4,\text{ }3,\text{ }2,\text{ }1 \right),
\quad{{\left( w_{4}^{3}\left( x \right) \right)}_{i}}=\left( 1,\text{ }3,\text{ }6,\text{ }10,\text{ }12,\text{ }12,\text{ }10,\text{ }6,\text{ }3,\text{ }1 \right),$$
$${{\left( w_{4}^{4}\left( x \right) \right)}_{i}}=\left( 1,4,10,20,31,40,44,40,31,20,10,4,1 \right);$$
$${{W}_{(1,\text{ }4)}}=\left( 4 \right),\qquad{{W}_{(2,\text{ }4)}}=\left( \begin{matrix}
   10 & 6  \\
   6 & 10  \\
\end{matrix} \right) ,\qquad{{W}_{(3,\text{ }4)}}=\left( \begin{matrix}
   20 & 10 & 4  \\
   40 & 44 & 40  \\
   4 & 10 & 20  \\
\end{matrix} \right).$$
$${{\left( w_{2}^{5}\left( x \right) \right)}_{i}}=\left( 1,\text{ }5,\text{ }10,\text{ }10,\text{ }5,\text{ }1 \right),$$
$${{\left( w_{3}^{5}\left( x \right) \right)}_{i}}=\left( 1,\text{ }5,\text{ }15,\text{ }30,\text{ }45,\text{ }51,\text{ }45,\text{ }30,\text{ }15,\text{ }5,\text{ }1 \right);$$
$${{W}_{(4,\text{ }2)}}=\left( \begin{matrix}
   5 & 1 & 0 & 0  \\
   10 & 10 & 5 & 1  \\
   1 & 5 & 10 & 10  \\
   0 & 0 & 1 & 5  \\
\end{matrix} \right) , \qquad{{W}_{(4,\text{ }3)}}=\left( \begin{matrix}
   15 & 5 & 1 & 0  \\
   51 & 45 & 30 & 15  \\
   15 & 30 & 45 & 51  \\
   0 & 1 & 5 & 15  \\
\end{matrix} \right).$$

Note the identities
$${{W}_{\left( n,\text{ }m \right)}}{{\tilde{A}}_{n}}\left( x \right)={{m}^{n}}{{\tilde{A}}_{n}}\left( x \right),$$
$${{W}_{\left( n,\text{ }m \right)}}{{\tilde{I}}_{n}}={{\tilde{I}}_{n}}{{W}_{\left( n,\text{ }m \right)}},  \qquad{{W}_{\left( n,\text{ }m \right)}}{{W}_{\left( n,\text{ }p \right)}}={{W}_{\left( n,\text{ }mp \right)}},$$
$${{W}_{\left( n,\text{ }m \right)}}{{\left( 1-x \right)}^{p}}{{c}_{n-p}}\left( x \right)={{\left( 1-x \right)}^{p}}{{W}_{\left( n-p,\text{ }m \right)}}{{c}_{n-p}}\left( x \right),  \qquad p<n,$$ 
where ${{c}_{n-p}}\left( x \right)$ is the polynomial of degree $<n-p$,  or
$$\left( {{\left( 1-x \right)}^{-p}},x \right){{W}_{\left( n,\text{ }m \right)}}\left( {{\left( 1-x \right)}^{p}},x \right){{I}_{n-p}}={{W}_{\left( n-p,\text{ }m \right)}}.$$
For example,
$$\left( \begin{matrix}
   4 & 1 & 0  \\
   4 & 6 & 4  \\
   0 & 1 & 4  \\
\end{matrix} \right)\left( \begin{matrix}
   1  \\
   4  \\
   1  \\
\end{matrix} \right)=8\left( \begin{matrix}
   1  \\
   4  \\
   1  \\
\end{matrix} \right),  \quad\left( \begin{matrix}
   10 & 4 & 1  \\
   16 & 19 & 16  \\
   1 & 4 & 10  \\
\end{matrix} \right)\left( \begin{matrix}
   1  \\
   4  \\
   1  \\
\end{matrix} \right)=27\left( \begin{matrix}
   1  \\
   4  \\
   1  \\
\end{matrix} \right),$$
$$\left( \begin{matrix}
   4 & 1 & 0  \\
   4 & 6 & 4  \\
   0 & 1 & 4  \\
\end{matrix} \right)\left( \begin{matrix}
   4 & 1 & 0  \\
   4 & 6 & 4  \\
   0 & 1 & 4  \\
\end{matrix} \right)=\left( \begin{matrix}
   20 & 10 & 4  \\
   40 & 44 & 40  \\
   4 & 10 & 20  \\
\end{matrix} \right),$$
$$\left( \begin{matrix}
   1 & 0 & 0  \\
   1 & 1 & 0  \\
   1 & 1 & 1  \\
\end{matrix} \right)\left( \begin{matrix}
   4 & 1 & 0  \\
   4 & 6 & 4  \\
   0 & 1 & 4  \\
\end{matrix} \right)\left( \begin{matrix}
   1 & 0  \\
   -1 & 1  \\
   0 & -1  \\
\end{matrix} \right)=\left( \begin{matrix}
   3 & 1  \\
   1 & 3  \\
\end{matrix} \right),$$
$$\left( \begin{matrix}
   1 & 0 & 0  \\
   2 & 1 & 0  \\
   3 & 2 & 1  \\
\end{matrix} \right)\left( \begin{matrix}
   4 & 1 & 0  \\
   4 & 6 & 4  \\
   0 & 1 & 4  \\
\end{matrix} \right)\left( \begin{matrix}
   1  \\
   -2  \\
   1  \\
\end{matrix} \right)=\left( \begin{matrix}
   1 & 0  \\
   1 & 1  \\
\end{matrix} \right)\left( \begin{matrix}
   3 & 1  \\
   1 & 3  \\
\end{matrix} \right)\left( \begin{matrix}
   1  \\
   -1  \\
\end{matrix} \right)=\left( 2 \right).$$

Since
$$\left( 1,a\left( x \right)-1 \right)\left( 1,{{\left( 1+x \right)}^{m}}-1 \right)=\left( 1,{{a}^{m}}\left( x \right)-1 \right),$$
matrix ${{W}_{\left( n,m \right)}}$ can also be represented in the form
$${{W}_{\left( n,m \right)}}=V_{n}^{-1}{{\left( \frac{{{\left( 1+x \right)}^{m}}-1}{x},{{\left( 1+x \right)}^{m}}-1 \right)}^{T}}{{V}_{n}}.$$
For example,
  $$\left( \begin{matrix}
   4 & 1 & 0  \\
   4 & 6 & 4  \\
   0 & 1 & 4  \\
\end{matrix} \right)=\left( \begin{matrix}
   1 & 0 & 0  \\
   -2 & 1 & 0  \\
   1 & -1 & 1  \\
\end{matrix} \right)\left( \begin{matrix}
   2 & 1 & 0  \\
   0 & 4 & 4  \\
   0 & 0 & 8  \\
\end{matrix} \right)\left( \begin{matrix}
   1 & 0 & 0  \\
   2 & 1 & 0  \\
   1 & 1 & 1  \\
\end{matrix} \right).$$
{\bfseries Theorem 5.} \emph{ Sum of the elements of each column of the matrix ${{W}_{\left( n,m \right)}}$ is ${{m}^{n}}$. }
\\
{\bfseries Proof.} According to the Theorem 2, $\alpha _{n}^{\left( m \right)}\left( 1 \right)={{\left( {{a}_{1}}m \right)}^{n}}$.
\\
{\bfseries Example 5.}
$$a\left( x \right)={{\left( 1-x \right)}^{-1}},  \qquad{{\tilde{\alpha }}_{n}}\left( x \right)=1,  \qquad\tilde{\alpha }_{n}^{\left( m \right)}\left( x \right)=[\uparrow ,0]{{W}_{\left( n,\text{ }m \right)}}.$$
In particular,
$$\tilde{\alpha }_{n}^{\left( 2 \right)}\left( x \right)=[\uparrow ,0]{{W}_{\left( n,\text{ }2 \right)}}={{r}_{n}}\prod\limits_{m=1}^{\left[ {n}/{2}\; \right]}{\left( x+\text{t}{{\text{g}}^{2}}\frac{m}{n+1}\pi  \right)},$$  
where ${{r}_{n}}=1$ for even $n$, ${{r}_{n}}=n+1$ for odd $n$, 
$$x\tilde{\alpha }_{n}^{\left( 2 \right)}\left( {{x}^{2}} \right)=\frac{{{\left( 1+x \right)}^{n+1}}-{{\left( 1-x \right)}^{n+1}}}{2}.$$
This corresponds to the case $a\left( x \right)={{\left( 1-2x+{{x}^{2}} \right)}^{-1}}$ in Example 3.
\section{GEP and generalized Lagrange series}

It follows from the Lagrange series expansion  for arbitrary formal power series $b\left( x \right)$ and $a\left( x \right)$, ${{a}_{0}}=1$:
$$\frac{b\left( x \right)}{1-x{{\left( \log a\left( x \right) \right)}^{\prime }}}=\sum\limits_{n=0}^{\infty }{\frac{{{x}^{n}}}{{{a}^{n}}\left( x \right)}}\left[ {{x}^{n}} \right]b\left( x \right){{a}^{n}}\left( x \right)$$
that each formal power series $a\left( x \right)$, ${{a}_{0}}=1$, is associated by means of the transform
$${{a}^{\varphi }}\left( x \right)=\sum\limits_{n=0}^{\infty }{\frac{{{x}^{n}}}{{{a}^{\beta n}}\left( x \right)}\left[ {{x}^{n}} \right]}\left( 1-x\beta {{\left( \log a\left( x \right) \right)}^{\prime }} \right){{a}^{\varphi +\beta n}}\left( x \right)$$
with the set of series$_{\left( \beta  \right)}a\left( x \right)$, $_{\left( 0 \right)}a\left( x \right)=a\left( x \right)$,   such that
$${}_{\left( \beta  \right)}a\left( x{{a}^{-\beta }}\left( x \right) \right)=a\left( x \right),   \qquad a\left( x{}_{\left( \beta  \right)}{{a}^{\beta }}\left( x \right) \right)={}_{\left( \beta  \right)}a\left( x \right),$$
$$\left[ {{x}^{n}} \right]{}_{\left( \beta  \right)}{{a}^{\varphi }}\left( x \right)=\left[ {{x}^{n}} \right]\left( 1-x\beta \frac{{a}'\left( x \right)}{a\left( x \right)} \right){{a}^{\varphi +\beta n}}\left( x \right)=\frac{\varphi }{\varphi +\beta n}\left[ {{x}^{n}} \right]{{a}^{\varphi +\beta n}}\left( x \right),$$
$$\left[ {{x}^{n}} \right]\left( 1+x\beta \frac{_{\left( \beta  \right)}{a}'\left( x \right)}{_{\left( \beta  \right)}a\left( x \right)} \right){}_{\left( \beta  \right)}{{a}^{\varphi }}\left( x \right)=\frac{\varphi +\beta n}{\varphi }\left[ {{x}^{n}} \right]{}_{\left( \beta  \right)}{{a}^{\varphi }}\left( x \right)=\left[ {{x}^{n}} \right]a_{{}}^{\varphi +\beta n}\left( x \right).$$
Series $_{\left( \beta  \right)}a\left( x \right)$ for integer $\beta $, denoted by ${{S}_{\beta }}\left( x \right)$, were introduced in [9]. In [26] these series, called generalized Lagrange series, are associated with the following construction. Table whose $k$th row, $k=0$, $\pm 1$, $\pm 2$, … , corresponds to the series 
$${{a}^{\beta k}}\left( x \right), \qquad{{a}_{0}}=1, \qquad\beta >0,$$
will be denoted by ${{\left\{ {{a}^{\beta }}\left( x \right) \right\}}_{0}}$. Table whose $k$th row is the $k$th ascending diagonal of the table ${{\left\{ {{a}^{\beta }}\left( x \right) \right\}}_{0}}$ will be denoted by ${{\left\{ {{a}^{\beta }}\left( x \right) \right\}}_{1}}$. Table whose $k$th row is the $k$th ascending diagonal of the table ${{\left\{ {{a}^{\beta }}\left( x \right) \right\}}_{1}}$ will be denoted by ${{\left\{ {{a}^{\beta }}\left( x \right) \right\}}_{2}}$; etc. For example, ${{\left\{ 1+x \right\}}_{0}}$, ${{\left\{ 1+x \right\}}_{1}}$, ${{\left\{ 1+x \right\}}_{2}}$:,  
$$\begin{matrix}
   \vdots   \\
   3  \\
   2  \\
   1  \\
   k=0  \\
   -1  \\
   -2  \\
   -3  \\
   \vdots   \\
\end{matrix}\left( \begin{matrix}
   \vdots  & \vdots  & \vdots  & \vdots  &   \\
   1 & 3 & 3 & 1 & \cdots   \\
   1 & 2 & 1 & 0 & \cdots   \\
   1 & 1 & 0 & 0 & \cdots   \\
   1 & 0 & 0 & 0 & \cdots   \\
   1 & -1 & 1 & -1 & \cdots   \\
   1 & -2 & 3 & -4 & \cdots   \\
   1 & -3 & 6 & -10 & \cdots   \\
   \vdots  & \vdots  & \vdots  & \vdots  & \ddots   \\
\end{matrix} \right),$$
  $$\begin{matrix}
   \vdots   \\
   3  \\
   2  \\
   1  \\
   k=0  \\
   -1  \\
   -2  \\
   -3  \\
   \vdots   \\
\end{matrix}\left( \begin{matrix}
   \vdots  & \vdots  & \vdots  & \vdots  &   \\
   1 & 4 & 10 & 20 & \cdots   \\
   1 & 3 & 6 & 10 & \cdots   \\
   1 & 2 & 3 & 4 & \cdots   \\
   1 & 1 & 1 & 1 & \cdots   \\
   1 & 0 & 0 & 0 & \cdots   \\
   1 & -1 & 0 & 0 & \cdots   \\
   1 & -2 & 1 & 0 & \cdots   \\
   \vdots  & \vdots  & \vdots  & \vdots  & \ddots   \\
\end{matrix} \right),  \qquad\begin{matrix}
   \vdots   \\
   3  \\
   2  \\
   1  \\
   k=0  \\
   -1  \\
   -2  \\
   -3  \\
   \vdots   \\
\end{matrix}\left( \begin{matrix}
   \vdots  & \vdots  & \vdots  & \vdots  &   \\
   1 & 5 & 21 & 84 & \cdots   \\
   1 & 4 & 15 & 56 & \cdots   \\
   1 & 3 & 10 & 35 & \cdots   \\
   1 & 2 & 6 & 20 & \cdots   \\
   1 & 1 & 3 & 10 & \cdots   \\
   1 & 0 & 1 & 4 & \cdots   \\
   1 & -1 & 0 & 1 & \cdots   \\
   \vdots  & \vdots  & \vdots  & \vdots  & \ddots   \\
\end{matrix} \right).$$
Table whose $k$th row is the $k$th descending diagonal of the table ${{\left\{ {{a}^{\beta }}\left( x \right) \right\}}_{0}}$ will be denoted by ${{\left\{ {{a}^{\beta }}\left( x \right) \right\}}_{-1}}$. Table whose $k$th row is the $k$th descending diagonal of the table ${{\left\{ {{a}^{\beta }}\left( x \right) \right\}}_{-1}}$ will be denoted by ${{\left\{ {{a}^{\beta }}\left( x \right) \right\}}_{-2}}$; etc. For example, ${{\left\{ 1+x \right\}}_{-1}}$, ${{\left\{ 1+x \right\}}_{-2}}$:
$$\begin{matrix}
   \vdots   \\
   3  \\
   2  \\
   1  \\
   k=0  \\
   -1  \\
   -2  \\
   -3  \\
   \vdots   \\
\end{matrix}\left( \begin{matrix}
   \vdots  & \vdots  & \vdots  & \vdots  &   \\
   1 & 2 & 0 & 0 & \cdots   \\
   1 & 1 & 0 & -1 & \cdots   \\
   1 & 0 & 1 & -4 & \cdots   \\
   1 & -1 & 3 & -10 & \cdots   \\
   1 & -2 & 6 & -20 & \cdots   \\
   1 & -3 & 10 & -35 & \cdots   \\
   1 & -4 & 15 & -56 & \cdots   \\
   \vdots  & \vdots  & \vdots  & \vdots  & \ddots   \\
\end{matrix} \right) ,   \qquad\begin{matrix}
   \vdots   \\
   3  \\
   2  \\
   1  \\
   k=0  \\
   -1  \\
   -2  \\
   -3  \\
   \vdots   \\
\end{matrix}\left( \begin{matrix}
   \vdots  & \vdots  & \vdots  & \vdots  &   \\
   1 & 1 & 1 & -10 & \cdots   \\
   1 & 0 & 3 & -20 & \cdots   \\
   1 & -1 & 6 & -35 & \cdots   \\
   1 & -2 & 10 & -56 & \cdots   \\
   1 & -3 & 15 & -84 & \cdots   \\
   1 & -4 & 21 & -120 & \cdots   \\
   1 & -5 & 28 & -165 & \cdots   \\
   \vdots  & \vdots  & \vdots  & \vdots  & \ddots   \\
\end{matrix} \right).$$
It turns out that the $k$th row of the table ${{\left\{ {{a}^{\beta }}\left( x \right) \right\}}_{v}}$ corresponds to the series 
$$\left( 1+xv\beta {{\left( \log {}_{\left( v\beta  \right)}a\left( x \right) \right)}^{\prime }} \right){}_{\left( v\beta  \right)}{{a}^{\beta k}}\left( x \right),$$
which follows from the identity
$$\left[ {{x}^{n}} \right]{{a}^{\beta \left( k+vn \right)}}\left( x \right)=\left[ {{x}^{n}} \right]\left( 1+xv\beta {{\left( \log {}_{\left( v\beta  \right)}a\left( x \right) \right)}^{\prime }} \right){}_{\left( v\beta  \right)}{{a}^{\beta k}}\left( x \right).$$

Denote
$$\left[ n,\to  \right]\left( 1,{}_{\left( \beta  \right)}a\left( x \right) \right)=\frac{_{\left( \beta  \right)}{{\alpha }_{n}}\left( x \right)}{{{\left( 1-x \right)}^{n+1}}},  \qquad\frac{1}{x}{}_{\left( \beta  \right)}{{\alpha }_{n}}\left( x \right)={}_{\left( \beta  \right)}{{\tilde{\alpha }}_{n}}\left( x \right),$$
$$\left[ n,\to  \right]{{\left( 1,\log {}_{\left( \beta  \right)}a\left( x \right) \right)}_{{{e}^{x}}}}={}_{\left( \beta  \right)}{{u}_{n}}\left( x \right),  \qquad\frac{1}{x}{}_{\left( \beta  \right)}{{u}_{n}}\left( x \right)={}_{\left( \beta  \right)}{{\tilde{u}}_{n}}\left( x \right).$$
Then
$$_{\left( \beta  \right)}{{a}^{\varphi }}\left( x \right)=\sum\limits_{n=0}^{\infty }{\frac{\varphi }{\varphi +n\beta }}\frac{{{u}_{n}}\left( \varphi +n\beta  \right)}{n!}{{x}^{n}},\quad
_{\left( \beta  \right)}{{u}_{n}}\left( x \right)=x{{\left( x+n\beta  \right)}^{-1}}{{u}_{n}}\left( x+n\beta  \right),$$
$${{E}^{n\beta }}{{\tilde{u}}_{n}}\left( x \right)={{\tilde{u}}_{n}}\left( x+n\beta  \right)={}_{\left( \beta  \right)}{{\tilde{u}}_{n}}\left( x \right),\quad
{{U}_{n}}{{E}^{n\beta }}U_{n}^{-1}{{\tilde{\alpha }}_{n}}\left( x \right)={}_{\left( \beta  \right)}{{\tilde{\alpha }}_{n}}\left( x \right).$$

Denote
$${{U}_{n}}{{E}^{n}}U_{n}^{-1}={{A}_{n}}.$$
Since
$$\left( 1,-x \right){{E}^{n}}\left( 1,-x \right)={{E}^{-n}},   \qquad{{U}_{n}}\left( 1,-x \right)U_{n}^{-1}={{\left( -1 \right)}^{n+1}}{{\tilde{I}}_{n}},$$
then
$${{\tilde{I}}_{n}}{{A}_{n}}{{\tilde{I}}_{n}}=A_{n}^{-1}.$$
For example,
$${{A}_{2}}=\left( \begin{matrix}
  2 & 1  \\
   -1 & 0  \\
\end{matrix} \right),   \quad{{A}_{3}}=\left( \begin{matrix}
 5 & {5}/{2}\; & 1  \\
   -6 & -2 & 0  \\
  2 & {1}/{2}\; & 0  \\
\end{matrix} \right),    
\quad{{A}_{4}}=\left( \begin{matrix}
14 & 7 & 3 & 1  \\
   -28 & -{35}/{3}\; & -{10}/{3}\; & 0  \\
  20 & {22}/{3}\; & {5}/{3}\; & 0  \\
   -5 & -{5}/{3}\; & -{1}/{3}\; & 0  \\
\end{matrix} \right);$$
$$A_{2}^{-1}=\left( \begin{matrix}
   0 & -1  \\
   1 & 2  \\
\end{matrix} \right),   \quad A_{3}^{-1}=\left( \begin{matrix}
   0 & {1}/{2}\; & 2  \\
   0 & -2 & -6  \\
   1 & {5}/{2}\; & 5  \\
\end{matrix} \right),   \quad A_{4}^{-1}=\left( \begin{matrix}
   0 & -{1}/{3}\; & -{5}/{3}\; & -5  \\
   0 & {5}/{3}\; & {22}/{3}\; & 20  \\
   0 & -{10}/{3}\; & -{35}/{3}\; & -28  \\
   1 & 3 & 7 & 14  \\
\end{matrix} \right).$$

Denote
$$A_{n}^{\beta }={{U}_{n}}{{E}^{n\beta }}U_{n}^{-1}.$$
For example,
$$A_{2}^{{1}/{2}\;}=\frac{1}{2}\left( \begin{matrix}
   3 & 1  \\
   -1 & 1  \\
\end{matrix} \right),  \quad A_{3}^{{1}/{2}\;}=\frac{1}{8}\left( \begin{matrix}
   21 & 7 & 1  \\
   -18 & 2 & 6  \\
   5 & -1 & 1  \\
\end{matrix} \right),  \quad A_{4}^{{1}/{2}\;}=\frac{1}{6}\left( \begin{matrix}
   30 & 10 & 2 & 0  \\
   -45 & -5 & 5 & 3  \\
   27 & 1 & -1 & 3  \\
   -6 & 0 & 0 & 0  \\
\end{matrix} \right).$$
{\bfseries Theorem 6.} \emph{ Sum of the elements of each column of the matrix $A_{n}^{\beta }$ is $1$.}
\\
{\bfseries Proof.} Since $\left[ x \right]{}_{\left( \beta  \right)}a\left( x \right)=\left[ x \right]a\left( x \right)={{a}_{1}}$, from the Theorem 2 it follows that ${}_{\left( \beta  \right)}{{\tilde{\alpha }}_{n}}\left( 1 \right)={{\tilde{\alpha }}_{n}}\left( 1 \right)$.
\\
{\bfseries Remark 3 (corollary of Remark 1).} If ${{c}_{n-m}}\left( x \right)$, $m<n$, is the polynomial of degree $<n-m$, then
$$A_{n}^{\beta }{{\left( 1-x \right)}^{m}}{{c}_{n-m}}\left( x \right)={{\left( 1-x \right)}^{m}}A_{n-m}^{{n\beta }/{\left( n-m \right)}\;}{{c}_{n-m}}\left( x \right),$$
or
$$\left( {{\left( 1-x \right)}^{-m}},x \right)A_{n}^{\beta }\left( {{\left( 1-x \right)}^{m}},x \right){{I}_{n-m}}=A_{n-m}^{{n\beta }/{\left( n-m \right)}\;}.$$

Denote
$$\log {{A}_{n}}={{U}_{n}}nDU_{n}^{-1}.$$
where $D$ is the matrix of the differential operator. Since
$${{E}^{n\beta }}=\sum\limits_{m=0}^{\infty }{\frac{{{\left( n\beta D \right)}^{m}}}{m!}},  \qquad nD=\log {{E}^{n}},$$
then
$$A_{n}^{\beta }=\sum\limits_{m=0}^{n-1}{\frac{{{\beta }^{m}}}{m!}}{{\left( \log {{A}_{n}} \right)}^{m}}.$$
For example,
$$A_{2}^{\beta }={{I}_{2}}+\beta \left( \begin{matrix}
  1 & 1  \\
   -1 & -1  \\
\end{matrix} \right),$$
$$A_{3}^{\beta }={{I}_{3}}+\beta \frac{1}{2}\left( \begin{matrix}
  5 & 2 & -1  \\
   -6 & 0 & 6  \\
   1 & -2 & -5  \\
\end{matrix} \right)+\frac{{{\beta }^{2}}}{2!}3\left( \begin{matrix}
   1 & 1 & 1  \\
   -2 & -2 & -2  \\
   1 & 1 & 1  \\
\end{matrix} \right),$$

$$A_{4}^{\beta }={{I}_{4}}+\beta \frac{1}{3}\left( \begin{matrix}
   13 & 3 & -1 & 1  \\
   -18 & 4 & 8 & -6  \\
   6 & -8 & -4 & 18  \\
   -1 & 1 & -3 & -13  \\
\end{matrix} \right)+\frac{{{\beta }^{2}}}{2!}\frac{4}{3}\left( \begin{matrix}
   9 & 5 & 1 & -3  \\
   -21 & -9 & 3 & 15  \\
   15 & 3 & -9 & -21  \\
   -3 & 1 & 5 & 9  \\
\end{matrix} \right)+$$
$$+\frac{{{\beta }^{3}}}{3!}16\left( \begin{matrix}
  1 & 1 & 1 & 1  \\
   -3 & -3 & -3 & -3  \\
   3 & 3 & 3 & 3  \\
   -1 & -1 & -1 & -1  \\
\end{matrix} \right),$$
where
$$\left[ \uparrow ,p \right]{{\left( \log {{A}_{n}} \right)}^{n-1}}={{n}^{n-2}}{{\left( 1-x \right)}^{n-1}}.$$
{\bfseries Example 6.} If $a\left( x \right)=1+x$, then $_{\left( \beta  \right)}a\left( x \right)$ is the generalized binomial series:
$$_{\left( \beta  \right)}{{a}^{\varphi }}\left( x \right)=\sum\limits_{n=0}^{\infty }{\frac{\varphi }{\varphi +n\beta }}\left( \begin{matrix}
   \varphi +n\beta   \\
   n  \\
\end{matrix} \right){{x}^{n}};$$
$$_{\left( 1 \right)}a\left( x \right)=\frac{1}{1-x},  \qquad_{\left( 2 \right)}a\left( x \right)=\frac{1-{{\left( 1-4x \right)}^{{1}/{2}\;}}}{2x},$$
$$_{\left( -1 \right)}a\left( x \right)=\frac{1+{{\left( 1+4x \right)}^{{1}/{2}\;}}}{2},     \qquad_{\left( {1}/{2}\; \right)}a\left( x \right)={{\left( \frac{x}{2}+{{\left( 1+\frac{{{x}^{2}}}{4} \right)}^{{1}/{2}\;}} \right)}^{2}}.$$
Since ${{\tilde{\alpha }}_{n}}\left( x \right)={{x}^{n-1}}$, then 
$$[\uparrow ,n-1]A_{n}^{\beta }={}_{\left( \beta  \right)}{{\tilde{\alpha }}_{n}}\left( x \right).$$
In particular, as follows from Example 2,
$$[\uparrow ,2n-1]A_{2n}^{{1}/{2}\;}=\frac{1}{2}\left( 1+x \right){{x}^{n-1}}.$$

We can come to the transformation $A_{n}^{\beta }$ in a different way, which leads to a simpler method of constructing the matrix $A_{n}^{\beta }$. We introduce the matrices $\tilde{D}=D\left( x,x \right)$, ${{\tilde{D}}^{-1}}$: 
$$\tilde{D}=\left( \begin{matrix}
   1 & 0 & 0 & \cdots   \\
   0 & 2 & 0 & \cdots   \\
   0 & 0 & 3 & \cdots   \\
   \vdots  & \vdots  & \vdots  & \ddots   \\
\end{matrix} \right),   \qquad{{\tilde{D}}^{-1}}=\left( \begin{matrix}
   1 & 0 & 0 & \cdots   \\
   0 & \frac{1}{2} & 0 & \cdots   \\
   0 & 0 & \frac{1}{3} & \cdots   \\
   \vdots  & \vdots  & \vdots  & \ddots   \\
\end{matrix} \right).$$
Let
 $$\frac{{{{\tilde{\alpha }}}_{n}}\left( x \right)}{{{\left( 1-x \right)}^{n+1}}}=\sum\limits_{m=0}^{\infty }{{{b}_{m}}}{{x}^{m}},  \qquad{{b}_{m}}=\left[ {{x}^{n}} \right]{{a}^{m+1}}\left( x \right).$$
Since
$$\left[ {{x}^{n}} \right]{}_{\left( 1 \right)}{{a}^{m+1}}\left( x \right)=\frac{m+1}{m+1+n}\left[ {{x}^{n}} \right]{{a}^{m+1+n}}\left( x \right),$$
then
$$\frac{_{\left( 1 \right)}{{{\tilde{\alpha }}}_{n}}\left( x \right)}{{{\left( 1-x \right)}^{n+1}}}=\sum\limits_{m=0}^{\infty }{\frac{m+1}{m+1+n}{{b}_{m+n}}{{x}^{m}}}.$$
Since
$$\tilde{D}{{\left( {{x}^{n}},x \right)}^{T}}{{\tilde{D}}^{-1}}\sum\limits_{m=0}^{\infty }{{{b}_{m}}}{{x}^{m}}=\tilde{D}{{\left( {{x}^{n}},x \right)}^{T}}\sum\limits_{m=0}^{\infty }{\frac{1}{m+1}{{b}_{m}}{{x}^{m}}}=$$
$$=\tilde{D}\sum\limits_{m=0}^{\infty }{\frac{1}{m+1+n}{{b}_{m+n}}{{x}^{m}}}=\sum\limits_{m=0}^{\infty }{\frac{m+1}{m+1+n}}{{b}_{m+n}}{{x}^{m}},$$
then
$$_{\left( 1 \right)}{{\tilde{\alpha }}_{n}}\left( x \right)=\left( {{\left( 1-x \right)}^{n+1}},x \right)\tilde{D}{{\left( {{x}^{n}},x \right)}^{T}}{{\tilde{D}}^{-1}}\left( {{\left( 1-x \right)}^{-n-1}},x \right){{\tilde{\alpha }}_{n}}\left( x \right),$$
$${{A}_{n}}=\left( {{\left( 1-x \right)}^{n+1}},x \right)\tilde{D}{{\left( {{x}^{n}},x \right)}^{T}}{{\tilde{D}}^{-1}}\left( {{\left( 1-x \right)}^{-n-1}},x \right){{I}_{n}}.$$
We use the identity  
$$D\left( x,x \right)\left( g\left( x \right),xg\left( x \right) \right)=D\left( 1,xg\left( x \right) \right)\left( x,x \right)=\left( {{\left( xg\left( x \right) \right)}^{\prime }},xg\left( x \right) \right)D\left( x,x \right),$$
or
$$\tilde{D}\left( g\left( x \right),xg\left( x \right) \right)=\left( {{\left( xg\left( x \right) \right)}^{\prime }},xg\left( x \right) \right)\tilde{D},$$
applied to the Pascal matrix, $g\left( x \right)={{\left( 1-x \right)}^{-1}}$ :
$$\tilde{D}P=\left( {{\left( 1-x \right)}^{-1}},x \right)P\tilde{D},$$
$$\left( 1-x,x \right)\tilde{D}=P\tilde{D}{{P}^{-1}},  \qquad{{\tilde{D}}^{-1}}\left( {{\left( 1-x \right)}^{-1}},x \right)=P{{\tilde{D}}^{-1}}{{P}^{-1}}.$$
Since 
$$E\left( {{x}^{n}},x \right){{E}^{-1}}=\left( {{\left( 1+x \right)}^{n}},x \right), \qquad{{P}^{-1}}{{\left( {{x}^{n}},x \right)}^{T}}P={{\left( {{\left( 1+x \right)}^{n}},x \right)}^{T}},$$
$$\left( {{\left( 1-x \right)}^{n}},x \right)P{{I}_{n}}=V_{n}^{-1},  \qquad{{P}^{-1}}\left( {{\left( 1-x \right)}^{-n}},x \right){{I}_{n}}=\left( {{\left( 1+x \right)}^{n}},x \right){{P}^{-1}}{{I}_{n}}={{V}_{n}},$$
we have:
$$A_{n}^{\beta }=\left( {{\left( 1-x \right)}^{n}},x \right)P\tilde{D}{{\left( {{\left( 1+x \right)}^{n\beta }},x \right)}^{T}}{{\tilde{D}}^{-1}}{{P}^{-1}}\left( {{\left( 1-x \right)}^{-n}},x \right){{I}_{n}}=$$
$$=V_{n}^{-1}\tilde{D}{{\left( {{\left( 1+x \right)}^{n\beta }},x \right)}^{T}}{{\tilde{D}}^{-1}}{{V}_{n}}.$$
For example,
$${{A}_{2}}=\left( \begin{matrix}
   1 & 0  \\
   -1 & 1  \\
\end{matrix} \right)\left( \begin{matrix}
   1 & 0  \\
   0 & 2  \\
\end{matrix} \right)\left( \begin{matrix}
   1 & 2  \\
   0 & 1  \\
\end{matrix} \right)\left( \begin{matrix}
   1 & 0  \\
   0 & \frac{1}{2}  \\
\end{matrix} \right)\left( \begin{matrix}
   1 & 0  \\
   1 & 1  \\
\end{matrix} \right),$$
$${{A}_{3}}=\left( \begin{matrix}
   1 & 0 & 0  \\
   -2 & 1 & 0  \\
   1 & -1 & 1  \\
\end{matrix} \right)\left( \begin{matrix}
   1 & 0 & 0  \\
   0 & 2 & 0  \\
   0 & 0 & 3  \\
\end{matrix} \right)\left( \begin{matrix}
   1 & 3 & 3  \\
   0 & 1 & 3  \\
   0 & 0 & 1  \\
\end{matrix} \right)\left( \begin{matrix}
   1 & 0 & 0  \\
   0 & \frac{1}{2} & 0  \\
   0 & 0 & \frac{1}{3}  \\
\end{matrix} \right)\left( \begin{matrix}
   1 & 0 & 0  \\
   2 & 1 & 0  \\
   1 & 1 & 1  \\
\end{matrix} \right),$$
$${{A}_{4}}=\left( \begin{matrix}
   1 & 0 & 0 & 0  \\
   -3 & 1 & 0 & 0  \\
   3 & -2 & 1 & 0  \\
   -1 & 1 & -1 & 1  \\
\end{matrix} \right)\left( \begin{matrix}
   1 & 0 & 0 & 0  \\
   0 & 2 & 0 & 0  \\
   0 & 0 & 3 & 0  \\
   0 & 0 & 0 & 4  \\
\end{matrix} \right)\left( \begin{matrix}
   1 & 4 & 6 & 4  \\
   0 & 1 & 4 & 6  \\
   0 & 0 & 1 & 4  \\
   0 & 0 & 0 & 1  \\
\end{matrix} \right)\left( \begin{matrix}
   1 & 0 & 0 & 0  \\
   0 & \frac{1}{2} & 0 & 0  \\
   0 & 0 & \frac{1}{3} & 0  \\
   0 & 0 & 0 & \frac{1}{4}  \\
\end{matrix} \right)\left( \begin{matrix}
   1 & 0 & 0 & 0  \\
   3 & 1 & 0 & 0  \\
   3 & 2 & 1 & 0  \\
   1 & 1 & 1 & 1  \\
\end{matrix} \right).$$

Denote
$$\left[ n,\to  \right]\left( 1,{}_{\left( \beta  \right)}a\left( x \right)-1 \right)={}_{\left( \beta  \right)}{{v}_{n}}\left( x \right), \qquad\frac{1}{x}{}_{\left( \beta  \right)}{{v}_{n}}\left( x \right)={}_{\left( \beta  \right)}{{\tilde{v}}_{n}}\left( x \right).$$
Then
$${{V}_{n}}{{U}_{n}}{{E}^{n\beta }}U_{n}^{-1}V_{n}^{-1}{{\tilde{v}}_{n}}\left( x \right)=\tilde{D}{{\left( {{\left( 1+x \right)}^{n\beta }},x \right)}^{T}}{{\tilde{D}}^{-1}}{{\tilde{v}}_{n}}\left( x \right)={}_{\left( \beta  \right)}{{\tilde{v}}_{n}}\left( x \right).$$
\\
{\bfseries Example 7.} This example was considered in [13] for a particular cases of the generalized binomial series. We will consider it from a more general point of view, using the transformations $\tilde{D}{{\left( {{\left( 1+x \right)}^{n\beta }},x \right)}^{T}}{{\tilde{D}}^{-1}}{{I}_{n}}$, $V_{n}^{-1}$. Let $a\left( x \right)=1+x$. Then $_{\left( \beta  \right)}a\left( x \right)-1={{x}_{\left( \beta  \right)}}{{a}^{\beta }}\left( x \right)$,  ${{\tilde{v}}_{n}}\left( x \right)={{x}^{n-1}}$,
$$_{\left( \beta  \right)}{{\tilde{v}}_{n}}\left( x \right)=\tilde{D}{{\left( {{\left( 1+x \right)}^{n\beta }},x \right)}^{T}}{{\tilde{D}}^{-1}}{{x}^{n-1}}=\sum\limits_{m=0}^{n-1}{\frac{m+1}{n}\left( \begin{matrix}
   n\beta   \\
   n-m-1  \\
\end{matrix} \right)}{{x}^{m}},$$
so that
$$\left( _{\left( \beta  \right)}{{a}^{\beta }}\left( x \right),{{x}_{\left( \beta  \right)}}{{a}^{\beta }}\left( x \right) \right)={{\tilde{D}}^{-1}}A\tilde{D},$$ 
where
$$\left[ n,\to  \right]A=\left[ n,\to  \right]\left( {{\left( 1+x \right)}^{\left( n+1 \right)\beta }},x \right).$$
For example, when $\beta =2$,
$${{\tilde{D}}^{-1}}\left( \begin{matrix}
   1 & 0 & 0 & 0 & 0 & \cdots   \\
   4 & 1 & 0 & 0 & 0 & \cdots   \\
   15 & 6 & 1 & 0 & 0 & \cdots   \\
   56 & 28 & 8 & 1 & 0 & \cdots   \\
   210 & 120 & 45 & 10 & 1 & \cdots   \\
   \vdots  & \vdots  & \vdots  & \vdots  & \vdots  & \ddots   \\
\end{matrix} \right)\tilde{D}=\left( \begin{matrix}
   1 & 0 & 0 & 0 & 0 & \cdots   \\
   2 & 1 & 0 & 0 & 0 & \cdots   \\
   5 & 4 & 1 & 0 & 0 & \cdots   \\
   14 & 14 & 6 & 1 & 0 & \cdots   \\
   42 & 48 & 27 & 8 & 1 & \cdots   \\
   \vdots  & \vdots  & \vdots  & \vdots  & \vdots  & \ddots   \\
\end{matrix} \right).$$
Since
$$\left[ m,\to  \right]V_{n}^{-1}=\sum\limits_{i=0}^{m}{\left( \begin{matrix}
   m-n  \\
   m-i  \\
\end{matrix} \right){{x}^{i}}}=\sum\limits_{i=0}^{m}{{{\left( -1 \right)}^{m-i}}}\left( \begin{matrix}
   n-i-1  \\
   m-i  \\
\end{matrix} \right){{x}^{i}},$$
then
$$\left[ {{x}^{m}} \right]{}_{\left( \beta  \right)}{{\tilde{\alpha }}_{n}}\left( x \right)=\left[ {{x}^{m}} \right]V{{_{n}^{-1}}_{\left( \beta  \right)}}{{\tilde{v}}_{n}}\left( x \right)=$$
$$=\sum\limits_{i=0}^{m}{{{\left( -1 \right)}^{m-i}}}\left( \begin{matrix}
   n-i-1  \\
   m-i  \\
\end{matrix} \right)\frac{\left( i+1 \right)}{n}\left( \begin{matrix}
   n\beta   \\
   n-i-1  \\
\end{matrix} \right)\frac{\left( n\beta -n+m+1 \right)!}{\left( n\beta -n+m+1 \right)!}=$$
$$=\frac{1}{n}\left( \begin{matrix}
   n\beta   \\
   n-m-1  \\
\end{matrix} \right)\sum\limits_{i=0}^{m}{{{\left( -1 \right)}^{m-i}}\left( i+1 \right)}\left( \begin{matrix}
   n\beta -n+m+1  \\
   m-i  \\
\end{matrix} \right)=$$ 
$$=\frac{1}{n}\left( \begin{matrix}
   n\beta   \\
   n-m-1  \\
\end{matrix} \right){{\left( -1 \right)}^{m}}\left( \begin{matrix}
   n\beta -n+m-1  \\
   m  \\
\end{matrix} \right)=\frac{1}{n}\left( \begin{matrix}
   n\beta   \\
   n-m-1  \\
\end{matrix} \right)\left( \begin{matrix}
   n\left( 1-\beta  \right)  \\
   m  \\
\end{matrix} \right).$$
Thus,
$$_{\left( \beta  \right)}{{\alpha }_{n}}\left( x \right)=\frac{1}{n}\sum\limits_{m=1}^{n}{\left( \begin{matrix}
   n\left( 1-\beta  \right)  \\
   m-1  \\
\end{matrix} \right)\left( \begin{matrix}
   n\beta   \\
   n-m  \\
\end{matrix} \right){{x}^{m}}}.$$
Note that
$$_{\left( 1-\beta  \right)}a\left( x \right)={}_{\left( \beta  \right)}{{a}^{-1}}\left( -x \right),  \qquad_{\left( 1-\beta  \right)}{{\alpha }_{n}}\left( x \right)=x{{\hat{I}}_{n}}{}_{\left( \beta  \right)}{{\alpha }_{n}}\left( x \right).$$

E-mail: {evgeniy\symbol{"5F}burlachenko@list.ru}

\begin{thebibliography}{99}


\bibitem{11} L. Shapiro, S. Getu, W. Woan, L. Woodson, The Riordan group, Discrete Appl. Math. 34 (1991) 229-339.
\bibitem{12} R. Sprugnoli, Riordan arrays and combinatorial sums, Discrete Math.132 (1994) 267-290.
\bibitem{13} D. Merlini, D. G. Rogers, R. Sprugnoli and M. C. Verri, On some alternative characterizations of Riordan arrays, Can. J. Math. 49 (1997) 301-320.
\bibitem{14} W, Wang, T. Wang, Generalized Riordan arrays, Discrete Math. 308 (2008) 6466-6500.
\bibitem{15}P. Barry, A study of integer sequences, Riordan arrays, Pascal-like arrays and Hankel transforms, University College Cork, 2009.
\bibitem{16} V. E. Hoggatt, Convolution triangles for generalized Fibonacci numbers, Fibonacci Quart. 8 (1970) 158-171.
\bibitem{17} V. E. Hoggatt and Marjorie Bicknell, Convolution triangles, Fibonacci Quart. 10 (1972) 599-609.
\bibitem{18} V. E. Hoggatt, and G. E. Bergum, Generalized convolution arrays, Fibonacci Quart. 13 (1975)  193-198.
\bibitem{19} V. E. Hoggatt and Paul S. Bruckman, H-convolution transform, Fibonacci Quart. 13 (1975) 357-368.
\bibitem{20} V. E. Hoggatt and Marjorie Bicknell, Pascal, Catalan, and general sequence convolution arrays in a matrix, Fibonacci Quart. 14 (1976) 135-142.
\bibitem{21} V. E. Hoggatt and Marjorie Bicknell, Sequences of matrix inverses from Pascal, Catalan, and related convolution arrays, Fibonacci Quart. 14 (1976) 224-232.
\bibitem{22} V. E. Hoggatt and Marjorie Bicknell, Catalan and related sequences arising from inverses of Pascal's triangle matrices, Fibonacci Quart. 14 (1976) 395-404.
\bibitem{23} V. E. Hoggatt and Marjorie Bicknell-Johnson, Numerator polynomial coefficient arrays for Catalan and related sequence convolution triangles, Fibonacci Quart. 15 (1977) 30-34.
\bibitem{24} G. E. Bergum and V. E. Hoggatt, An application of the characteristic of the generalized Fibonacci sequence,  Fibonacci Quart. 15 (1977) 215-220.
\bibitem{25} V. E. Hoggatt and Marjorie Bicknell-Johnson, Properties of generating functions of a convolution array, Fibonacci Quart. 16 (1978) 289-295.
\bibitem{26} V. E. Hoggatt and Marjorie Bicknell-Johnson, Convolution arrays for Jacobsthal and Fibonacci polynomials, 16 (1978) 385-402.
\bibitem{27} Peter Bala, Notes on generalized Riordan arrays, https://oeis.org/A260492/a260492.pdf
\bibitem{28} M. Koutras, Eulerian nambers associated with sequences of polynomials, Fibonacci Quart. 32 (1994) 44-57.
\bibitem{29} T. K. Petersen, Eulerian Numbers, Birkhauser, 2015.
\bibitem{30} J. M. Holte, Carries, combinatorics, and an amazing matrix, Am. Math. Mon. 104 (1997) 138-149.
\bibitem{31} F. Brenti, V. Welker, The Veronese construction for formal power series and graded algebras,  Adv. in Appl. Math. 42 (2009) 545-556.
\bibitem{32} P. Diaconis, J. Fulman, Foulkes characters, Eulerian idempotents, and an amazing matrix, J. Algebr. Comb. 36 (2012) 425-440.
\bibitem{33} J.-C. Novelli, J.-Y. Thibon, Noncommutative symmetric functions and an amazing matrix, Adv. in Appl. Math. 48 (2012) 528-534.
\bibitem{34} L. Carlitz and V. E. Hoggatt, Generalized Eulerian numbers and polynomials, Fibonacci Quart. 16 (1978) 138-146
\bibitem{35} E. Burlachenko, Algebra of formal power series, isomorphic to the algebra of formal Dirichlet series, arXiv:1702.01071.
\bibitem{36} E. V. Burlachenko, Riordan arrays and generalized Lagrange series, Mathematical Notes, 100 (2016) 531-539.



\end{thebibliography}
\end{document}